\documentclass{article}

\usepackage[english]{babel}
\usepackage{authblk}
\usepackage{pdflscape}
\usepackage{ws-rotating}

\usepackage[letterpaper,top=2cm,bottom=2cm,left=3cm,right=3cm,marginparwidth=1.75cm]{geometry}

\usepackage{amsmath, tabu}
\usepackage{amsfonts}
\usepackage{array}
\usepackage{graphicx}
\usepackage{calrsfs}
\usepackage[table]{xcolor}
\usepackage{booktabs}

\usepackage{hyperref}
\hypersetup{verbose}
\usepackage{enumitem}
\setlist[enumerate]{itemsep=0mm}
\usepackage{multirow}
\hypersetup{colorlinks=false,allbordercolors=blue,pdfborderstyle={/S/U/W 1}}
\usepackage{tikz}

\providecommand{\keywords}[2]{\textit{Keywords---} #1}

\begin{document}

\markboth{\textit{M. Martin et al.}}{XPPLORE: Import, visualize, and analyze XPPAUT data in MATLAB}

\title{XPPLORE: \\Import, visualize, and analyze XPPAUT data in MATLAB}

\author[]{Matteo Martin\footnote{Author for correspondence - matteo.martin.2@phd.unipd.it}\footnote{These authors contributed equally to this work}}
\author[2]{Anna Kishida Thomas$^\dag$}
\author[2]{G. Bard Ermentrout}

\affil[1]{Department of Information Engineering, University of Padova}
\affil[2]{Department of Mathematics, University of Pittsburgh}

\maketitle

\begin{abstract}
\noindent The analysis of ordinary differential equation (ODE) dynamical systems, particularly in applied disciplines such as mathematical biology and neuroscience, often requires flexible computational workflows tailored to model-specific questions. XPPAUT is a widely used tool combining numerical integration and continuation methods. Various XPPAUT toolboxes have emerged to customize analyses, however, they typically rely on summary `.dat' files and cannot parse the more informative `.auto' files, which contain detailed continuation data, e.g. periodic orbits and boundary value problem solutions. We present XPPLORE, a user-friendly and structured MATLAB toolbox overcoming this limitation through the handling of `.auto' files. This free software enables post-processing of continuation results, facilitates analyses such as manifold reconstruction and averaging, and it supports the creation of high-quality visualizations suitable for scientific publications. This paper introduces the core data structures of XPPLORE and demonstrates the software's exploration capabilities, highlighting its value as a customizable and accessible extension for researchers working with ODE-based dynamical systems.
\end{abstract}

\hspace{0.35cm}\keywords{XPPAUT; AUTO; numerical continuation; ODE; bifurcation diagram}

\newpage
\tableofcontents
\thispagestyle{empty}

%\begin{multicols}{2}

\newpage
\section{Introduction}

In the field of dynamical systems, new questions and models arise daily, motivating novel analysis methods. Most of these innovations require computational and analytical exploration with customized workflows and \textit{ad-hoc} tools. In applied disciplines such as mathematical biology and neuroscience, investigating the dynamical mechanisms in a new model is crucial for understanding the biological system, the experimental results and for stating new predictions~\cite{Martin2024, Nechyporenko2024,JohnJP2024,JohnChaos2024}. A large component of this work deals with continuous deterministic models defined by systems of ordinary differential equations (ODEs). The most widely used tools to perform analyses on ODE-based dynamical systems are AUTO-07p~\cite{AUTO07p}, XPPAUT~\cite{XPPAUT}, MATCONT~\cite{MatCont} and PYDSTool~\cite{PyDSTool}.

Among these, XPPAUT combines XPP with the powerful continuation engine AUTO. This allows XPPAUT to efficiently compute one- and two-parameter bifurcation diagrams (1P-/2P-BDs) through numerical continuation. For more details on continuation techniques, the reader is referred to~\cite{Kuznetsov2023}. The results of these calculations can be exported in specific file formats (`.auto', `.dat') containing different levels of detail. 
While the capabilities of XPPAUT are vast, researchers have developed MATLAB and Python packages to meet their individual project needs. These toolboxes facilitate custom parameter and initial condition sweeps, visualizations of bifurcation diagrams, simulation, and calculations of nullclines \cite{XPPy, PlotBD, Py_XPPCALL, PyXPP, xppMATLAB, xppToolbox}. Despite these extensions of functionality, they can only handle the simplest format of continuation calculations (`.dat'), limiting their scope.

The continuation calculations in `.auto' files are the most valuable outputs of XPPAUT. They contain many details on computed bifurcation diagrams, including the numerical approximation of limit cycles, solutions to boundary-value problems (BVPs), and the specific settings used for the computation. So far, there has been no simple way to fully utilize the information in these files. This gap motivated the development of our package XPPLORE (pronounced ``explore"), a structured and interactive MATLAB \cite{MATLAB} toolbox for handling the XPPAUT continuation, simulation, and calculation of nullclines. Its main advantage is the parsing of the `.auto' file, opening the door to the reconstruction of manifolds, averaging, and other post-processing analyses. Another advantage lies in the ease of use of this package, as it does not require any changes to the XPPAUT workflow, but simply leverages output files from the creation of bifurcation diagrams for post-processing analyses. In addition, the toolbox enables users to generate high-quality images for scientific publications, making it an invaluable resource for students, researchers, and others working with ODE-based dynamical systems.

This paper introduces the core data structures of XPPLORE and illustrates its capabilities through example demonstrations.The document is divided into four sections: The first section covers the technical details of XPPAUT. The second presents the main data structures implemented in XPPLORE, while the third shows some examples illustrating a selection of the toolbox's capabilities. We conclude with final remarks on the toolbox.

\newpage
\section{XPPAUT file types}

\noindent XPPAUT~\cite{XPPAUT} is a computational tool to simulate, analyze, and animate dynamical systems. It incorporates a version of AUTO~\cite{AUTO07p}, software created to perform numerical continuation. In this section, we present the
files that XPPLORE handles regarding the model, data, nullclines, and continuation. Table~\ref{Table:ExportXPPAUT} illustrates the procedures for exporting these output files from XPPAUT.

\begin{table}[ht]
    \centering
    \caption{Table illustrating how to export simulations, nullclines and numerical continuation calculations from XPPAUT. \label{Table:ExportXPPAUT}}
    {\tabcolsep13pt
    \begin{tabular}{cc}\\[-2pt]
        \toprule
        Type & Procedure \\[6pt]
        \hline
        Data & Data $\rightarrow$ Write \\
        Nullclines & Nullcline $\rightarrow$ Save \\
        Continuation & File $\rightarrow$ Save diagram\\
        \hline
    \end{tabular}}
\end{table}

We warn the reader that a `.dat' file is a generic file type exportable from XPPAUT containing different types of results (e.g., bifurcation diagrams, numerical integrations, and nullclines). However, XPPLORE can only parse the `.dat' files containing results of numerical integrations or calculations of nullclines. There are numerous other packages mentioned in the introduction that are able to use the .dat files exported for bifurcation diagrams.

\subsection{Model file (`.ode')} \label{Section:Modelfile}

This file contains an implementation of the dynamical system,
\begin{equation}\label{EQ:ODESYS}
    \dot{x} = f(x\,|\,\theta)
\end{equation}
where $x$ and $\theta$ represent the $n$- and $m$-dimensional real state and parameter vectors. The right-hand side is defined as $f: \mathbb{R}^n \times \mathbb{R}^m \rightarrow \mathbb{R}^n$ with $f=[f_1,f_2,...,f_n]^T$. The dot notation indicates the derivative with respect to time variable $t$.
In the `.ode' files, the parameters, variables, boundary conditions, functions and numerical/visualization settings are explicitly stated.

\subsection{Data file (`.dat')} \label{Section:Datfile}

This file contains the results associated with numerical integration of system (\ref{EQ:ODESYS}) for a finite interval of time. The data are organized in a table whose columns, proceeding from left to right, present the time, dynamical, and auxiliary variables of the model. The ordering of the variables within these three subgroups is consistent with the one used in the `.ode' file with auxiliary quantities always following the dynamical variables.

\subsection{Nullclines file (`.dat')} \label{Section:Nullclinesfile}

This file contains the numerical approximations of the $x_i$- and $x_j$-nullclines of system (\ref{EQ:ODESYS}). These are the zero level sets of the functions $f_{i/j}(x\,|\,\theta)$ in the ($x_i$, $x_j$) plane. % where the condition =0$ is satisfied.
For those systems with $n > 2$, during the calculation of these sets of points, XPPAUT treats the variables non-visualized in the current graphics view as parameters, fixed to the values specified in the initial conditions windows.

\subsection{Continuation file (`.auto')} \label{Section:Continuationfile}

The `.auto' files store continuation calculations. Numerical continuation is a powerful method to track fixed points (equilibria) and limit cycles of (\ref{EQ:ODESYS}) by changing the system parameters. More specifically, this technique can be used to compute bifurcation diagrams, which are used to investigate the invariant sets of system (\ref{EQ:ODESYS}), their stability, and how they qualitatively and quantitatively evolve in response to changes in parameters. XPPAUT handles up to two-parameter continuations. The calculated results are saved and exported using the command \verb|Save| \verb|diagram| (see Table~\ref{Table:ExportXPPAUT}). Most existing MATLAB and Python extension packages of XPPAUT handle the results exported through the commands \verb|Write| \verb|pts| and \verb|All|\ \verb|info|  \cite{PlotBD, xppToolbox}. These files significantly limit the post-processing analyses, as they only contain the minimum and maximum values of limit cycles or solutions to BVP, rather than containing full approximations of the system trajectories.

The `.auto' file is organized into three parts: settings, points, and special solutions. The first section stores the AUTO numerical and visualization settings specified when the results are saved. The second section contains all bifurcation diagrams/continuation points, in the order they were computed. Finally, the third section presents the approximations of the system's fixed points, limit cycles, and BVP solutions associated with so-called labeled points (more details in the notes below).

% TECHNIQUALITIES
Here are a few notes regarding the numerical continuation in XPPAUT. 
\begin{enumerate}[label=(\roman*)]
    \item Reset of storage: %We would like to note that 
    The storage of numerically approximated solutions %associated with limit cycles and BVPs in XPPAUT
    is cleared every time a new one-parameter continuation is started.
    \item Labeled points: These are: bifurcation, user-defined, and regular points. The bifurcation points are automatically identified during numerical continuation through embedded test functions. User-defined points satisfy one or more active conditions declared under \verb|Usr| \verb|period| menu. Finally, the regular points are those points computed every \verb|NPR| steps. The user can define custom conditions or modify the \verb|NPR| setting to control the set of points at which numerical approximated solutions are calculated.
    %how many numerical approximations of limit cycles or solutions to BVP will be computed.
    \item Multiple bifurcation diagrams: XPPAUT can handle multiple bifurcation diagrams at the same time. Moreover, it can save their results in a single `.auto' file. We warn the reader that, in accordance with point (i) in this list, the numerically approximated solutions exported will be related to the last 1P-BD or to the last couple of 1P- and 2P-BD.
    \item Hot parameters: These are the eight (or fewer) parameters, found under the \verb|Parameters| window of AUTO, exported in the `.auto' files.
\end{enumerate}

\section{Data Structures}

The implemented data structures can be classified as primary and secondary structures. The former are those created by directly parsing XPPAUT files, such as the model (M), the data (DATA), the nullclines (NC), and the AUTO repository (AR). The latter can be retrieved from primary structures and include the bifurcation diagram (BD), the eigenvalues (EIG) and the trajectories (TRJ). Their name recalls the type of data they store, and their structure design aims to simplify the readability and accessibility of the content.

Among the primary structures, M is used to create DATA and AR. For this reason, as we will show later in the text, one of the first operations to use our toolbox consists of creating the M structure. XPPLORE can only parse those results that can be exported from XPPAUT through the procedures presented in Table~\ref{Table:ExportXPPAUT}.

\subsection{Model} \label{Section:Model}
 
\begin{figure}[!ht]
    \centering
    \includegraphics{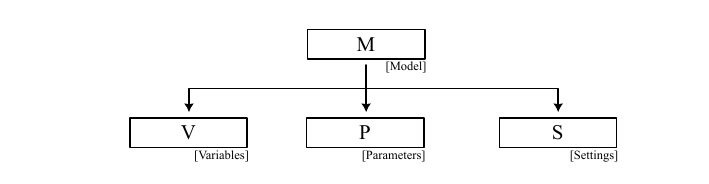}
    \caption{Overview of the M structure in XPPLORE. Its fields are \texttt{V}, \texttt{P}, and \texttt{S}, referring to variables, parameters, and settings, respectively.}
    \label{fig:Scheme:M}
    \vspace*{12pt}
\end{figure}

\noindent The M structure is created using the function \verb|Func_ReadModel()|, which parses the `.ode' file containing the dynamical system. This routine extracts the variables, parameters, and XPPAUT settings necessary to load the results of the computations. The M structure follows the organization depicted in Fig.~\ref{fig:Scheme:M}. In the following, the fields are illustrated.

\begin{enumerate}[label=(\roman*)]
    \item \verb|V|. A structure whose fields' names are the names of the temporal, dynamical, and auxiliary variables declared in the `.ode' file. Each field contains a character identifying the class of the variable (\verb|`T'|-temporal, \verb|`D'|-dynamical, \verb|`A'|-auxiliary).
    
    \item \verb|P|. A structure whose fields' names are the names of the model parameters. Each field contains the parameter value specified in the `.ode' file.

    \item \verb|S|. A structure storing the hard-coded XPPAUT settings. The fields of this structure match the names of the hard-coded settings, while the value of each field is the setting's hard-coded value.
\end{enumerate}

\subsection{Data} \label{Section:Simulations}

The DATA structure is created by parsing the `.dat' file using the function \verb|Func_ReadData()|. We remark that a `.dat' dataset file can only be used when exported through the procedure illustrated in Table~\ref{Table:ExportXPPAUT}. Therefore, a DATA structure can contain information regarding any kind of analysis requiring the simulation of the model, e.g. a plain integration or the results of Poincaré maps iterated over a finite parameter range. The function \verb|Func_ReadData()| produces a structure with fields whose names match those in M under \verb|V|. Hence, they are the names of the temporal, dynamical, and auxiliary variables. The content of each field is a vector storing the data in the corresponding column of the parsed file.

\subsection{Nullclines} \label{Section:Nullclines}

\noindent The `.dat' nullclines file can be loaded in MATLAB through \verb|Func_ReadNullclines()|. This function creates a structure whose fields are \verb|NCx| and \verb|NCy|, where \verb|x| and \verb|y| are replaced by the variables declared in the nullclines file name. That is, the `.dat' nullclines file name must have the following format: `[text]\_x\_y.dat', where `text' is any sequence of alpha-numerical characters, while x and y are the variables along the x and the y axis in the XPPAUT window, from which the nullclines are exported.

\subsection{AUTO repository} \label{Section:AutoRepo}

\begin{figure}[!ht]
    \centering
    \includegraphics{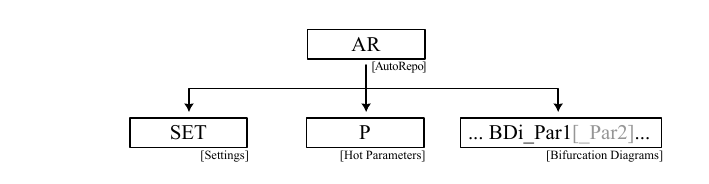}
    \caption{Visual representation of the AR structure and its fields \texttt{SET}, \texttt{P}, \texttt{BDi\_Par1[\_Par2]}, referring to continuation settings, model parameters and bifurcation diagram structures.} \label{Fig:Scheme:AR}
    \vspace*{12pt}
\end{figure}

The AUTO repository structure (AR) can be created by parsing an `.auto' file through the function \verb|Func_ReadAutoRepo()|. The output of this function is an AR structure whose fields are illustrated in the following paragraph and depicted in Fig.~\ref{Fig:Scheme:AR}. The AR structure mimics the organization of an `.auto' file.

\begin{enumerate}[label=(\roman*)]
\item \verb|SET|. Structures storing the settings used for the continuation with fields \verb|NPTS|, \verb|NUM|, and \verb|UZ|. 
\begin{enumerate}[label=(\alph*)]
    \item \verb|NPTS| stores the total number of points exported in the `.auto' file. 
    \item \verb|NUM| is a structure containing the numerical and visualization settings specified under the \verb|nUmerics| and \verb|Axes| menus of AUTO. 
    The fields of \verb|NUM| are: \verb|NTST|, \verb|NPR|, \verb|NMAX|, \verb|DS|, \verb|DSMIN|, \verb|DSMAX|, \verb|ParMIN|, \verb|ParMAX|, \verb|NormMIN|, \verb|NormMAX|, \verb|IAD|, \verb|MXBF|, \verb|IID|, \verb|ITMX|, \verb|ITNW|, \verb|NWTN|, \verb|IADS|, \verb|xmin|, \verb|ymin|, \verb|xmax|, and \verb|ymax|. For more information, please refer to the XPPAUT documentation~\cite{XPPAUT}.
    \item \verb|UZ| is a nested structure. It contains several fields as the number of user-defined conditions set under the \verb|Usr period| menu of AUTO. Their name has the form \verb|UZi| where \verb|i| $\in \mathbb{N}$, 1$\,\le\,$\verb|i|$\,\le\,$9. \verb|UZi| is a structure with a single field, where the field name corresponds to the parameter being evaluated in the condition, and its value represents the point at which the condition holds true.
\end{enumerate}

\item \verb|P|. A structure with fields \verb|Par1| to \verb|ParN|, with \verb|N| $\in \mathbb{N}$ and \verb|N| $\leq 8$. The field \verb|Pari| contains the i-th parameter declared under the \verb|Parameter| window of AUTO.

\item \verb|BDi_Par1[_Par2]|. This field is a structure corresponding to a bifurcation diagram stored in the `.auto' file. The name of this field presents an integer \verb|i| increasing across successive BDs. \verb|Par1| is replaced with the name of the main continuation parameter, while \verb|Par2| only exists in 2P-BDs, where it is replaced with the name of the secondary continuation parameter. 
\end{enumerate}

\subsubsection{Bifurcation diagram}\label{Section:BifDiag}

\begin{figure}[!ht]
    \centering
    \includegraphics{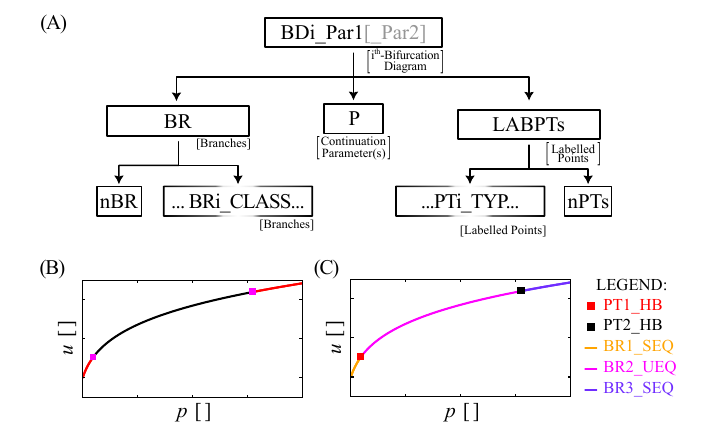}
    \caption{(A) Scheme summarizing the organization of the i-th bifurcation diagram computed using the continuation parameter(s) \texttt{Par1} (and \texttt{Par2}) saved in the overarching AR structure. \texttt{BR} and \texttt{LABPTs} are addressed in Fig.~\ref{Fig:Scheme:BRiPTi}. Panel (B) shows a generic 1P-BD (\texttt{BDi\_p}) stored in the overarching AR structure. In both panels, \textit{u} indicates a generic dynamical variable while \textit{p} is the main continuation parameter. Panel (B) shows an example of a BD computed through XPPAUT. The \textit{black}/\textit{red curves} represent unstable/stable fixed points, whereas the \textit{magenta dots} correspond to Hopf bifurcations (HBs). Panel (C) illustrates the classification and decomposition in different branches (\texttt{BRi\_CLASS}) and labeled points (\texttt{PTi\_TYP}) of the 1P-BD as indicated in the legend of the chart.
    }
    \label{Fig:Scheme:BD}
    \vspace*{12pt}
\end{figure}

\noindent The main fields of the AR structure are the bifurcation diagram structures \verb|BDi_Par1[_Par2]|. 
A bifurcation diagram (BD) is a set of equilibrium points and periodic orbits of an ODE system, obtained by changing one or more parameters, called bifurcation parameters, at a time.
XPPLORE, following this definition, considers bifurcation diagrams that use the same bifurcation parameter(s) but have different values in the remaining parameters as separate diagrams. In a bifurcation diagram, the information is organized into different branches and labeled points. Specifically, we define a branch within a BD as a collection of equilibrium points or periodic orbits, sequentially computed with the same stability properties. With this in mind, the organization of the BD structures is illustrated in Fig.~\ref{Fig:Scheme:BD}(A) and described in detail in the following paragraphs. Instead, Figs.~\ref{Fig:Scheme:BD}(B) and~\ref{Fig:Scheme:BD}(C) show an example of a 1P-BD and its decomposition into different branches and labeled points.

\begin{enumerate}[label=(\roman*)]

\item \verb|P|. Cell array containing the continuation parameters. In 1P-/2P-BDs, the array consists of one/two element(s).

\item \verb|BR|. A structure whose fields are the branches of the BD. The naming convention of the subfields follows \verb|Bi_CLASS|, where \verb|i| $\in \{1,...,\mathrm{nBR}\}$ is an incremental index and \verb|nBR| is the number of branches in the current BD. \verb|CLASS| indicates the branch type as specified in Table~\ref{Table_BR:Class}.

\begin{table}[ht]
    \centering
    \caption{Lists of values that CLASS can assume.\label{Table_BR:Class} EQ: equilibrium, LC: limit cycle, BVP: Boundary-value problem.}
    {\tabcolsep13pt\begin{tabular}{cc|cc}\\[-2pt]
        \toprule
        CLASS & Description & CLASS & Description \\[6pt]
        \hline
        SEQ & Stable EQ & SN & Saddle-node \\
        UEQ & Unstable EQ & SNPO & Saddle-node of periodic orbits\\
        SLC & Stable LC & HB & Andronov-Hopf \\
        ULC & Unstable LC & TR & Torus \\
        BVP & Solution of BVP & BP & Branch-point\\
        UZ & User-defined point & PD & Period-doubling\\
        \hline
    \end{tabular}}
\end{table}

Each branch \verb|Bi_CLASS| is a structure with multiple fields. Fig.~\ref{Fig:Scheme:BRiPTi} shows the organization of this structure. We indicate with $p$ the number of points in the considered branch. In the following, each field is explained in detail.

\begin{enumerate}[label=(\alph*)]

    \item \verb|BR|. An integer representing the 1-based index of the branch in the `.auto' file. 
    
    \item \verb|TPar|. An integer describing the type of continuation used for the current branch. If this value is \verb|0|, a single parameter is used. Otherwise, the branch is the result of a two-parameter continuation, and the value corresponds to the algorithm used from AUTO for the calculation. 
    
    \item \verb|TYP|. An integer representing the type of points stored in the current branch. Table~\ref{Table_BR:Fields} illustrates the different types.

    \begin{table}[ht]
        \centering
        \caption{Table illustrating the different types TYP of points that form branches in a bifurcation diagram.
        \label{Table_BR:Fields}}
        {\tabcolsep13pt\begin{tabular}{ccc|ccc}\\[-2pt]
            \toprule
            TPar & TYP & Description & TPar & TYP & Description \\[6pt]
            \hline
            0 & 1 & SEQ & 1 & 1 & SN   \\
            0 & 2 & UEQ & 2 & 2 & SNPO \\
            0 & 3 & SLC & 3 & 3 & HB   \\
            0 & 4 & ULC & 4 & 4 & TR   \\
            0 & 8 & BVP & 5 & 5 & BP   \\
            9 & 9 & UZ  & 6 & 6 & PD    \\
            \hline
        \end{tabular}}
    \end{table}

    \item \verb|IDX.| The absolute index of the points in the current branch relative to the full diagram.

    \item \verb|Par1...ParN.| The names of these fields are the names of the parameters under the field \verb|P| in \verb|AR|. Each of these fields contains a vector with $p$ entries. These vectors store information regarding the values of the hot parameters during the continuation.
    %with a number of entries as the number of points in the current branch.

    \item \verb|L2.| Array containing the $\mathcal{L}_2$-norm of each solution.

    \item \verb|T.| Array containing the periods of each solution. This field exists only in those branches associated with LCs or BVPs.

    \item \verb|ui/uU/uL/uA.| These are the initial/upper/lower/average values of a dynamical variable \verb|u| along each solution. Each branch contains these four fields for each dynamical variable involved in the system.

    \item \verb|EigR/EigI.| These fields consist of $n\times p$ -matrices.
    %, where $p$ is the number of points in the current branch and $n$ is the dimension of the system. 
    If the branch stores information regarding fixed points, the matrix \verb|EigR| and \verb|EigI| store the real and the imaginary components of the exponential of the eigenvalues. Instead, if the branch contains information on periodic orbits, \verb|EigR| and \verb|EigI| present the real and the imaginary components of the Floquet multipliers.

\end{enumerate}

\item \verb|LABPTs|. This field contains all labeled points detected during one- or two-parameter continuations. The structure has two subfields: \verb|nPT| and \verb|PTi_TYP|. The former stores the number of labeled points in the current BD, while the latter contains information about the i-th labeled point, whose type is \verb|TYP|. Table~\ref{Table_LABPTs:TYP} illustrates the values that \verb|TYP| can assume. 

\begin{table}
    \centering
    \caption{Table illustrating the TYP values that labeled points can assume. \label{Table_LABPTs:TYP}}
    {\tabcolsep13pt\begin{tabular}{cc|cc}\\[-2pt]
        \toprule
        TYP & Description & TYP & Description \\[6pt]
        \hline
        HB & HB & TR & TR \\
        SN & SN & EP & Endpoint \\
        PD & PD & UZ & User-defined point\\
        SNPO & SNPO& &\\
        \hline
    \end{tabular}}
\end{table}

Each labeled point \verb|PTi_TYP| is a structure containing the same fields of a branch, and two additional ones, as depicted in Fig.~\ref{Fig:Scheme:BRiPTi},

\begin{enumerate}[label=(\alph*)]
    \item \verb|LAB.| The XPPAUT label of the point.
    \item \verb|PO.| In the case of an LC or BVP, this field contains an approximation of the LC or of the solution of a BVP associated with the labeled point.
\end{enumerate}

\begin{figure}[!ht]
    \centering
    \includegraphics{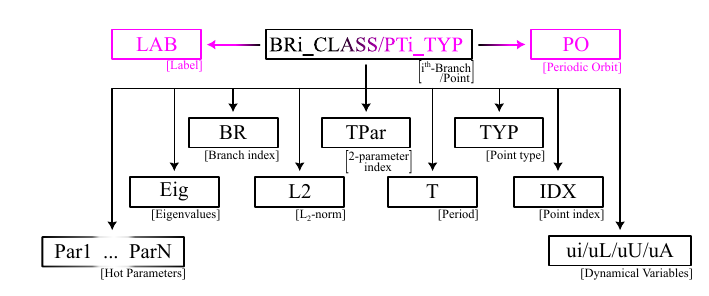}
    \caption{Organization of the i-th branch/labeled point structure within a bifurcation diagram structure. The fields in \textit{magenta} are contained only in labeled points associated with POs or BVPs.}
    \label{Fig:Scheme:BRiPTi}
    \vspace*{12pt}
\end{figure}

\end{enumerate}

\subsubsection{Eigenvalues} \label{Section:Eigenvalues}

\begin{figure}
    \centering
    \includegraphics{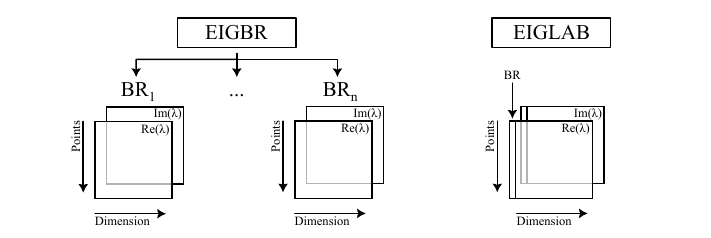}
    \caption{Overview of eigenvalue structures for branch points (\texttt{EIGBR}) and labeled points (\texttt{EIGLAB}). \texttt{EIGBR} contains one matrix for each branch within a BD, whereas \texttt{EIGLAB} contains a single matrix. $\lambda$: exponentials of the eigenvalues or Floquet multipliers of a point within a BD.}
    \label{Fig:Eigenvals}
    \vspace*{12pt}
\end{figure}

EIGBR and EIGLAB are secondary structures created using the function \verb|Func_GetEIG()|, which contain information regarding eigenvalues of the Jacobian matrix at equilibrium points and Floquet multipliers of periodic orbits. EIGBR presents this information sorted by branches within a BD, while EIGLAB contains these values only at labeled points. Even though the information of EIGLAB is contained in EIGBR, XPPLORE considers this subset separately because the specific labeled point can be hard to detect among all branch points. 
These structures are organized as shown in Fig.~\ref{Fig:Eigenvals}, and present the following details,

\begin{enumerate}[label=(\roman*)]
    \item \verb|EIGBR|. This structure contains $n\times p$-matrices, one for each branch of a considered BD. 
    Each matrix stores either the real or imaginary components of $\lambda$, that is, for an equilibrium point (limit cycle), the exponential of the eigenvalue (Floquet multiplier).
    \item \verb|EIGLAB|. It is a matrix containing the branch index and the eigenvalue or Floquet multiplier at bifurcation and user-defined points. Which branch this point is grouped with by XPPAUT solely depends on the numerical continuation and on which branch the point appears first. Note that, for example, in a system admitting a HB point, the bifurcation can belong to a stable or unstable branch of equilibria depending on where the continuation is started from.
\end{enumerate}

\subsubsection{Trajectories} \label{Section:Trajectories}

The final secondary structure in this package is TRJ, created by using the function \verb|Func_GetTRJ()|. This function extracts from a BD structure all the trajectories stored under the \verb|PO| field of labeled points. These trajectories can be LCs or solutions of BVPs. The fields of this structure are presented below.

\begin{figure}[ht!]
    \centering
    \includegraphics{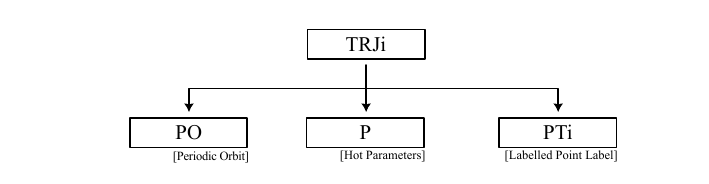}
    \caption{Illustration of the \texttt{TRJi} structure and its fields \texttt{PO}, \texttt{P}, and \texttt{PTi}.}
    \label{Fig:Scheme:TRJi}
    \vspace*{12pt}
\end{figure}

\begin{enumerate}[label=(\roman*)]
    \item \verb|TRJi.| This field stores a structure associated with the i-th extracted special trajectory. \verb|i| is an incremental index varying between 1 and \verb|nTRJ|. Its fields are \verb|PO|, \verb|P|, and \verb|PTi| as presented in Fig.~\ref{Fig:Scheme:TRJi} and described in detail below.
    \begin{enumerate}[label=(\alph*)]
        \item \verb|PO| is a structure with fields, normalized-by-period time \verb|t|, and the temporal evolution of the system's dynamical variables.
        \item \verb|P| is a parameter structure whose fields are the eight (or fewer) hot parameters, and the period \verb|T| of the solution. 
        \item \verb|PTi| is the name of the point in the \verb|LABPTs| field of a bifurcation diagram structure from where the solution has been extracted.
    \end{enumerate}
    \item \verb|nTRJ.| The number of \verb|TRJi| structures in \verb|TRJ|.
\end{enumerate}

%%%%%%%%%%%%%%%%%%%%%%%%%% XPPLORE-ation %%%%%%%%%%%%%%%%%%%%%%%%%%%%%%
\newpage
\section{XPPLORE-ation}\label{Section:XXPLOREation}

\subsection{Setup and preparation}
This toolbox requires working versions of both MATLAB \cite{MATLAB} and XPPAUT \cite{XPPAUT} as well as basic knowledge of these two programs. To use XPPLORE, the toolbox must be downloaded from \href{https://github.com/MatteoMartin/XPPLORE} {GitHub} and unzipped.\\

\noindent The reader is welcome to follow the examples in this section by referring to the provided `Demo' code folder. In all examples, the two operations of (1) importing XPPLORE functionalities and (2) making the toolbox available in the current MATLAB session are assumed to be fulfilled. These steps should be inserted at the beginning of an \textit{.m} or \textit{.mlx} file:

\begin{enumerate}[label=(\arabic*)]
    \item Environment initialization: Close all open figures, clear the workspace and clean the command window. The MATLAB commands are
    \vspace{-5pt}
    \begin{verbatim} clear all; close all; clc; \end{verbatim}
    \vspace{-5pt}
    \item Import XPPLORE: Using the \verb|addpath()| and \verb|genpath()| commands. \verb|PATH| is a variable containing the relative or absolute path to the XPPLORE folder.
    \vspace{-5pt}
    \begin{verbatim} addpath(genpath(PATH)); \end{verbatim}
\end{enumerate}

In the next sections, the one- or two-parameter bifurcation diagrams will be indicated using the notation $P_1$-BD or ($P_1$,$P_2$)-BD, where $P_{1/2}$ is the name of the main/secondary continuation parameter. 

\begin{figure}[ht!]
    \centering
    \includegraphics{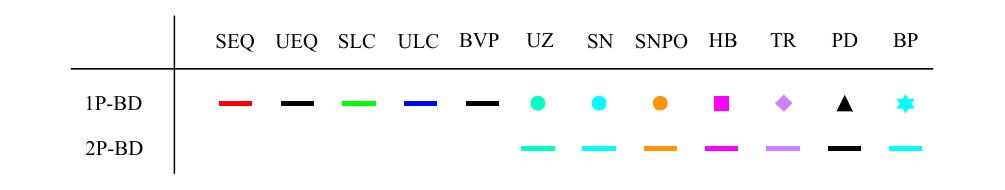}
    \caption{Legend of colors and shapes used in XPPLORE to represent EQ, LC, BVP, bifurcation points in 1P/2P-BDs, and nullclines in planar visualizations. These are default colors and shapes. For more information on how the user can change these, see the Appendix section.}
    \label{Fig:Scheme:LEGEND}
    \vspace*{12pt}
\end{figure}

The fixed points, bifurcations and nullclines in this section will follow the color scheme illustrated in Table~\ref{Fig:Scheme:LEGEND}.

% \begin{table}[ht]
%     \tbl{Lists of values that CLASS can assume.\label{Table_LEGEND:Class} EQ: equilibrium, LC: limit cycle, BVP: Boundary-value problem, SN: saddle-node bifurcation, SNPO: saddle-node of periodic orbits, HB: Andronov-Hopf bifurcation, TR: torus bifurcation, PD: period-doubling bifurcation.}
%     {\tabcolsep13pt\begin{tabular}{cc|cc}\\[-2pt]
%         \toprule
%         CLASS & Description & CLASS & Description \\[6pt]
%         \hline
%         \circlecolor{red}{} & Stable EQ & \circlecolor{red}{} & Saddle-node \\
%         \circlecolor{red}{} & Unstable EQ & \circlecolor{red}{} & SNPO\\
%         \circlecolor{red}{} & Stable LC & \circlecolor{red}{} & HB \\
%         \circlecolor{red}{} & Unstable LC & \circlecolor{red}{} & TR \\
%         \circlecolor{red}{} & Solution of BVP & \circlecolor{red}{} & PD\\
%         \circlecolor{red}{} & Fixed-period & & \\
%         \botrule
%     \end{tabular}}
% \end{table}

%%%%%%%%%%%%%%%%%%%%%%%%%% Model & Simulation %%%%%%%%%%%%%%%%%%%%%%%%%%

\subsection{Model, simulation and  nullclines}\label{Section:ModelSimNullclines}

\noindent\textbf{Demo}. The accompanying MATLAB file can be found in the folder \verb|DEMOs/DEMO1_ML_MDL|.

\vspace{5pt}
\noindent\textbf{Model.} To introduce model, simulation, and nullcline structures, we use a simplified version of the Morris-Lecar model~\cite{Morris1981}.
%~\cite{ermentrout1996reflected,morris1981voltage},
\begin{equation}
    \begin{aligned}
        c_m \dot{V} &= I_0 - I_L(V) - I_{Ca}(V,n) - I_{K}(V),\\
        \dot{n} &= \phi[n_{\infty}(V)-n] \tau_n^{-1}(V).
    \end{aligned}
\label{SYS:ML}
\end{equation}

\noindent The transmembrane voltage $V$ and gating variable $n$ are governed by the equations above, and the parameters are taken from~\cite{Morris1981}. With the aim of illustrating the utility of the nullclines to discuss the trajectory evolution in the phase space, we scale the parameter $\phi$ from 0.04 ms\textsuperscript{-1} by a factor of 10. The formulation for the fast persistent calcium current and delayed-rectifier potassium current is as follows,
\begin{equation}
    \begin{aligned}
        I_{Ca}(V,n) &= g_{Ca} \; n (V - E_{Ca}), \\
        I_{K}(V) &= g_K \; m_{\infty}(V) (V - E_K).
    \end{aligned}
    \label{SYS:EQ:MLCURRENTs}
\end{equation}
\noindent\textbf{Data.} The data considered in this example are:
\begin{enumerate}[label=(\roman*)]
    \item `.ode' file: File containing an implementation of system (\ref{SYS:ML}).
    \item `.dat' simulation file: Simulation of system (\ref{SYS:ML}) with $I_0=100\, \mu$A/cm\textsuperscript{2}. 
    %After numerically integrating the system in XPPAUT, export the data using the procedure illustrated in table~\ref{Table:ExportXPPAUT}.
    \item `.dat' nullcline file: $n$- and $V$-nullclines computed with $I_0=100\, \mu$A/cm\textsuperscript{2}. %Export this information through the \verb|Save| command under the \verb|Nullcline| menu. The name of the nullcline data file must obey the format presented in the documentation of the function \verb|Func_Read| \verb|Nullclines()|. 
\end{enumerate}

\noindent\textbf{Steps.} %The next few steps introduce how to import the model and few basic visualization operations in XPPLORE.
\begin{enumerate}[label=(\arabic*)]
%    \item Read the content of an `.ode' file
%    \item Read the content of a `.dat' simulation file
%    \item Read the content of a `.dat' nullclines file
%    \item Visualize a simulation
%    \item Visualize the nullclines
%\end{enumerate}

\item \textit{Read the content of an `.ode' file}.
To load any data, the user must create the model structure by parsing the `.ode' file with the function \verb|Func_ReadModel()|. The input argument to this function is the path (\verb|PATHODEFILE|) to the `.ode' file.
\begin{verbatim}
    M = Func_ReadModel(PATHODEFILE);
\end{verbatim}
This creates a model structure and saves it in the variable \verb|M|. For more details regarding its organization, please refer to Sec.~\ref{Section:Model}. The content of the variable \verb|M|, investigated through the MATLAB command window, is the following,
\begin{verbatim}
    M = 
        struct with fields:
        V: [1×1 struct]
        P: [1×1 struct]
        S: [1×1 struct]
\end{verbatim}
The model structure, as previously described in Sec.~\ref{Section:Model}, presents three fields: \verb|V|, \verb|P| and \verb|S|. \verb|V| stores the auxiliary and the dynamical variables implemented in the `.ode' file. \verb|P| contains the model parameters and their default values. \verb|S| is a structure of hard-coded numerical settings.

\item \textit{Read the content of a `.dat' simulation file.}
The simulation of a system can be loaded through the function \verb|Func_ReadData()|, which accepts as input the model structure (\verb|M|) and the simulation file name. Its output is a simulation structure (see Sec.~\ref{Section:Simulations}).
\begin{verbatim}
    D = Func_ReadData(M,`sim.dat');
\end{verbatim}
By typing \verb|D| in the MATLAB command window, the result is
\begin{verbatim}
    D = 
        struct with fields:
           t: [10001×1 double]
           v: [10001×1 double]
           n: [10001×1 double]
            :
        stim: [10001×1 double]
\end{verbatim}
Each of the fields contains a vector.

\item \textit{Read the content of a `.dat' nullclines file.} 
To read the nullclines exported through XPPAUT, the function \verb|Func_ReadNullclines()| can be used. The `.dat' nullcline file must respect the naming convention specified in Sec.~\ref{Section:Nullclines} and in the documentation of this function. For more information, we invite the reader to look at the appendix.
\begin{verbatim}
    NC = Func_ReadNullclines(`nc_n_v.dat');
\end{verbatim}

\item \textit{Visualize a simulation.}
For simulation visualization, XPPLORE uses the built-in MATLAB function \verb|plot()|, rather than implementing an ad-hoc one. The following line of code solves the task.
\begin{verbatim}
    fig = figure();
    plot(D.t,D.v,`Color',`b',`LineWidth',1.2)
\end{verbatim}
The result is presented in Fig.~\ref{Fig:DEMO1:DSWEB}(A).

\item \textit{Visualize the nullclines}. The nullclines can be visualized through the function \verb|Func_VisualizeNullclines()| which accepts a nullcline structure \verb|NC| as input.
\begin{verbatim}
    fig = figure();
    Func_VisualizeNullclines(NC);
\end{verbatim}
Fig.~\ref{Fig:DEMO1:DSWEB}(B) shows the nullclines of system (\ref{SYS:ML}) computed using $I_0$ set to 100 $\mu$A/cm\textsuperscript{2}, along with the simulated model trajectory obtained by solving the system using $V$ and $n$ I.C.s set to -27 mV and 0.14, respectively. By default, the function \verb|Func_VisualizeNullclines()| visualizes the nullclines by plotting over the x- and the y-axis the same dynamical variables visualized over the phase plane in the XPPAUT window used to compute the curves. However, the user can switch the axes of the variables displayed along the x- and y-axes by specifying the optional input \verb|`VAR'| to the visualization function.

\begin{figure}[ht]
    \centering
    \includegraphics{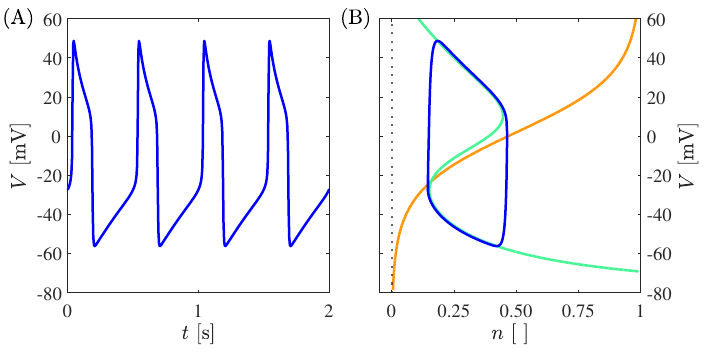}
    \caption{(A) Temporal evolution of system (\ref{SYS:ML}) for $I_0= 100\ \mu$A/cm\textsuperscript{2}. (B) The blue curve represents the system simulation. The orange and light blue curves are the $n$ and $V$ nullclines, respectively. The dotted black vertical line indicates $n=0$. The nullclines are computed by using $I_0$ set to $100\ \mu$A/cm\textsuperscript{2}.}
    \label{Fig:DEMO1:DSWEB}
    \vspace*{12pt}
\end{figure}
\end{enumerate}

%%%%%%%%%%%%%%%%%%%%%%%%%% BD & EVal %%%%%%%%%%%%%%%%%%%%%%%%

\subsection{Bifurcation diagrams} \label{Section:DEMOBD}

\noindent\textbf{Model}. This example considers the Hodgkin and Huxley model~\cite{HH1952}.
\begin{equation}
    \begin{aligned}
        c_m\dot{V} &= I_0 - I_{Na}(V,m,h) - I_K(V,n) - I_L(V), \\
        \dot{p} &= \alpha_p(V) (1-p) - \beta_p(V) p \:\:\: \mathrm{with} \:\:\:p \in\{m,n,h\}
    \end{aligned}
    \label{SYS:HH}
\end{equation}
where $V$ represents the membrane potential and $p$ indicates any gating variable. $\alpha_p(V),\ \beta_p(V)$ are the transition rates of $p$. $c_m$ is the membrane capacitance and $I_0$ represents the applied current. $I_{Na}$, $I_K$ and $I_L$ are the Na\textsuperscript{+}, K\textsuperscript{+} and leakage currents, respectively, defined as follows,
\begin{equation}
    \begin{aligned}
        I_{Na}(V,m,h) &= g_{Na} \; m^3 h (V - E_{Na}), \\
        I_{K}(V,n) &= g_{K} \; n^4 (V - E_K), \\
        I_{L}(V) &= g_L \; (V - E_L)
    \end{aligned}
    \label{SYS:EQ:CURRENTs}
\end{equation}
with $g$s and $E$s indicating the maximal conductance and Nernst potential, respectively.

% 1P-BD & Eigenvalues -------------------------------------------------

\subsubsection{One-parameter BD and eigenvalues} \label{Section:1PBD&Eig}

\noindent\textbf{Demo}. The accompanying MATLAB file can be found in the folder \verb|DEMOs/DEMO2_HH_1BD&EIG|.

\noindent\textbf{Data.} Data required for this section:
\begin{enumerate}[label=(\roman*)]
    \item `.ode' file containing the implementation of system (\ref{SYS:HH}).
    \item `.auto' file with continuation results of $I_0$-BD. %To obtain this file, compute the bifurcation diagram in XPPAUT and, in the AUTO window, select \verb|File| $
    %\rightarrow$ \verb|Save Diagram|.
\end{enumerate}

%\noindent\textbf{Overview.} 
\noindent\textbf{Steps.}
\begin{enumerate}[label=(\arabic*)]
%    \item Read the content of an `.auto' file
%    \item Visualize 1P-BD and its labeled points
%    \item Change the axes of a 1P-BD visualization
%    \item Visualize eigenvalues with 1P-BD
%    \item Retrieve eigenvalues for further use
%\end{enumerate}   

\item \textit{Preliminary step.} To perform the analyses presented in this section, the model structure \verb|M| has to be created. For more information, see Sec.~\ref{Section:ModelSimNullclines}.

\item \textit{Read the content of an `.auto' file.} The parsing of an `.auto' file can be done through the function \verb|Func_ReadAutoRepo()|, which creates an AR structure. It requires the model structure (\verb|M|, see Sec.~\ref{Section:ModelSimNullclines} for details) and the relative or absolute path to the `.auto' file (\verb|PATHAUTOFILE|). If the file is located in the MATLAB working directory, the name of the `.auto' file suffices as \verb|PATHTOFILE|.
\begin{verbatim}
    AR = Func_ReadAutoRepo(M,PATHAUTOFILE);
\end{verbatim}
The status of the function can be assessed in the MATLAB command window:
\begin{verbatim}
    AR:
    Parsing settings ...completed!
    Parsing hot parameters ...completed!
    Parsing branches & labeled points ...completed!
    Parsing Special Solutions...completed!

    Summary:
        1P-BD - Name: BD1_i0 - Main: i0
\end{verbatim}
The status output is organized into \verb|AR| and \verb|Summary|. The former shows the current task being executed and its status. While, in the latter, the BDs identified in the `.auto' file are summarized. In the previous example, the loaded file contains the $I_0$-BD only. The \verb|AR| structure investigated through the MATLAB command window is the following,
\begin{verbatim}
    AR = 
         struct with fields:
            SET: [1×1 struct]
              P: [1×1 struct]
         BD1_i0: [1×1 struct]
\end{verbatim}
As seen in Sec.~\ref{Section:AutoRepo}, the \verb|AR| structure has three or more fields, \verb|SET|, \verb|P|, and \verb|BDi_Par1[_Par2]|. %\verb|SET| contains the numerical and visualization settings used for the continuation, while \verb|P| contains the hot parameters. To explore the structure \verb|AR| in more detail, please refer to the associated \textit{.m} file. 

\item \textit{Visualize 1P-BD and its labeled points.} The bifurcation diagram and its labeled points can be visualized using \verb|Func_VisualizeDiagram()| and \verb|Func_VisualizeLabPoints()|, which require as input the model (\verb|M|) and a bifurcation diagram structure. The exact spelling of the bifurcation parameters (case-sensitive!) matches the one in the `.ode' file.
\begin{verbatim}
    fig = figure();
    Func_VisualizeDiagram(M,AR.BD1_i0)
    Func_VisualizeLabPoints(M,AR.BD1_i0)
    Func_FigStyle(fig)
\end{verbatim}
The boundaries and the labels of the axes in Fig.~\ref{Fig:Demo2:1PBD&Eig}(A), were adjusted through MATLAB customization functions. We warn the reader that, the function \verb|Func_VisualizeDiagram()| implements few optional fields as reported in the documentation of the function presented in the Appendix. Among these, \verb|TYPE| allows to control the visualization of the branch of periodics or the branch associated with solutions to BVPs. The default value of this hyperparameter is \verb|Standard| so that, for each solution, both the maximum and the minimum value of the y-axis variable are visualized. However, it can be changed to \verb|Lower|/\verb|Upper|/\verb|Initial|/\verb|Average| to visualize the minimum/maximum/initial/average value only of the specified variable along the solution.

\begin{figure}[ht]
    \centering
    \includegraphics{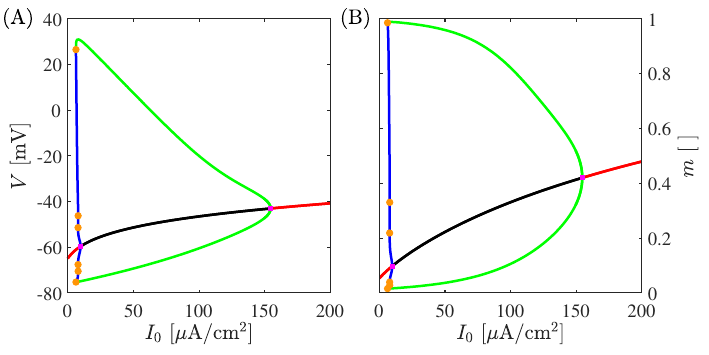}
    \caption{The continuous red and black lines represent stable and unstable fixed points, respectively. The green (stable) and the blue (unstable) lines are the maxima and minima of limit cycles. The orange dots are SNPOs while the magenta squares are HBs. (A) $I_0$-BD of system (\ref{SYS:HH}) visualized in the ($I_0$, $V$)-plane. (B) The same diagram of (A) but visualized over the ($I_0$, $m$)-plane. }\label{Fig:Demo2:1PBD&Eig}
    \vspace*{12pt}
\end{figure}

\item \textit{Change the axes of a 1P-BD visualization.} To change the visualized variable along the y-axis, we use the optional input argument \verb|VAR|, which must be specified through a name-value pair. Here, the dynamical variable $m$ is visualized against the bifurcation parameter $I_0$.
\begin{verbatim}
    fig = figure();
    Func_VisualizeDiagram(M,AR.BD1_i0,`VAR',{`i0',`m'})
    Func_VisualizeLabPoints(M,AR.BD1_i0,`VAR',{`i0',`m'})
\end{verbatim}
The result is presented in Fig.~\ref{Fig:Demo2:1PBD&Eig}(B).\\
\noindent Note that, only if a parameter is an `hot parameter' %- i.e. one of the eight parameters under the \verb|Parameters| tab in the AUTO window in XPPAUT when computing the BD - 
can it be used as a \verb|VAR| argument. See Sec.~\ref{Section:Continuationfile} for more information.

\item \textit{Visualize the period, the frequency or L\textsubscript{2}-norm.} The \verb|VAR| specifies which parameters, dynamical variables, period (\verb|T|), frequency (\verb|F|), or L\textsubscript{2}-norm (\verb|L2|) are to be visualized. In the case of visualizing the period/frequency of LC, or the L\textsubscript{2}-norm of EQs/LCs/BVPs, the \verb|VAR| field must have in the name-value pair \verb|{`i0',`T'}|/\verb|{`i0',`F'}|, or to \verb|{i0',`L2'}|, respectively. These capabilities are not demonstrated here.

\item \textit{Visualize exponential of eigenvalues and Floquet multipliers of a 1P-BD.} The function \\
\verb|Func_VisualizeEig()| shows how the real and the imaginary part of the exponential of eigenvalues or Floquet multipliers change when the bifurcation parameter is varied. The real and the imaginary components are visualized over the y- and the z-axis. The function additionally plots the cylinder with unitary radius, because the bifurcations of fixed points (in case of the exponential of eigenvalues) or periodic orbits occur at the unitary circle. The bifurcations are highlighted by coloring the section of the cylinder in correspondence with the parameter value where they occur. The color is consistent with the one used in BD.
%since the bifurcations of fixed point or periodic orbit occurs  of Floquet multipliers or exponentials of eigenvalues to easily identify the bifurcations.
%the planes at $\Re[\lambda]$ or $\Im[\lambda]$ equal to -1, 0 or 1, 
\begin{verbatim}
    Func_VisualizeEig(M,AR.BD1_i0)
\end{verbatim}
\noindent One can also specify a subset of branch indices within \verb|AR.BD1_i0.BR|. This optional argument must be specified as a name-value pair (e.g. \verb|`BRIND',{1,2,3,4}|). 
%Inspecting \verb|AR.BD1_i0.BR| reveals a total of 10 branches.
\begin{verbatim}
    Func_VisualizeEig(M,AR.BD1_i0,`BRIND',{1,2,3,4})
\end{verbatim}
\noindent In AUTO, labeled points are saved and assigned to one branch of the bifurcation diagram. Because of this, labeled points of deleted branches will be removed too.
\begin{verbatim}
    Func_VisualizeEig(M,AR.BD1_i0,`BRIND',{1,2,4})
\end{verbatim}
\noindent The function \verb|Func_VisualizeEig()| accepts as input the optional parameter \verb|VAR| that has to be specified through a name-value declaration. This field takes in input a cell vector whose entries specify the combination of components visualized along the x, y, and eventually z-axis of the chart. For more information on the variables accepted as input to this field, we ask the reader to refer to the appendix section.
%\noindent Lastly, we can visualize the imaginary components of the eigenvalues by specifying the \verb|TYPE| argument.
%\begin{verbatim}
%    Func_VisualizeEig(M,AR.BD1_i0,`TYPE',`Im')
%\end{verbatim}

\begin{figure}[ht]
    \centering
    \includegraphics{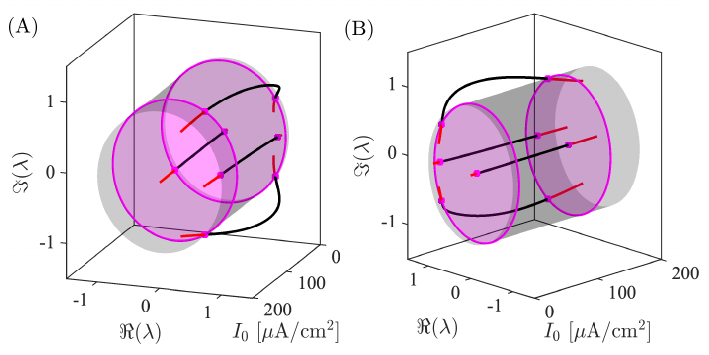}
    \caption{Visualization of the eigenvalues of all the fixed point branches of a $I_0$-BD. (A) and (B) visualize the same set of eigenvalues. The number of branches visualized corresponds to the number of eigenvalues computed for each fixed point. Two branches of eigenvalues cross the bifurcation cylinder at the HBs (magenta squares). The sections of the cylinder in correspondence of the HBs are highlighted in magenta.}\label{Fig:Demo2:Eig3d}
    \vspace*{12pt}
\end{figure}

\item \textit{Retrieving eigenvalues for further use.} Another way to make use of the eigenvalues of a 1P-BD is the function \verb|Func_GetEig()|. This function returns two objects, the first (second) is related to eigenvalues along branches (labeled points). More information regarding these structures can be found in Sec.~\ref{Section:Eigenvalues}. The example below retrieves the eigenvalues of the branches/labeled points into \verb|EIGBR|/\verb|EIGLAB| within \verb|AR.BD1_i0|, respectively.
\begin{verbatim}
    [EIGBR, EIGLAB] = Func_GetEig(AR.BD1_i0);
\end{verbatim}
\noindent For example, the below code accesses the 5th point within the 2nd branch of the bifurcation diagram, the real parts ("1") of all 4 (":") eigenvalues.
\begin{verbatim}
    EIGBR{2}(5,:,1) = [0.0095  0.7579  0.7579  0.8861]
\end{verbatim}
\noindent This code accesses, within the 3rd labeled point, the branch index and the imaginary ("2") components of all 4 eigenvalues.\\
\begin{verbatim}
    EIGLAB(3,:,2) = [3  0  0.8738  -0.8738  0]
\end{verbatim}

\end{enumerate}

% 1 & 2 Par BD -------------------------------------------------

\subsubsection{One- and two-parameter BDs}

\noindent\textbf{Demo}. The accompanying MATLAB file can be found in the folder \verb|DEMOs/DEMO3_HH_1P2PBDs|.

\noindent \textbf{Data.} For the analyses presented here, the following files are needed:
\begin{enumerate}[label=(\roman*)]
    \item \textit{.ode} file: File containing an implementation of system (\ref{SYS:HH}).
    \item \textit{.auto} file: File containing the $I_0$- and the ($I_0$, $g_K$)-BDs of system (\ref{SYS:HH}). In XPPAUT, first compute the $I_0$-BD, and then continue the most important bifurcations along the $g_K$ direction by creating the ($I_0$, $g_K$)-BD. After these two operations, save the \textit{.auto} file via \verb+File+ $
    \rightarrow$ \verb|Save| \verb|diagram|.
\end{enumerate}

\noindent \textbf{Steps.}
%\begin{itemize}[nolistsep]
%    \item[--] Read the content of an \textit{.auto} file
%    \item[--] Visualize a 2P-BD
%\end{itemize}

\begin{enumerate}[label=(\arabic*)]

\item \textit{Preliminary step.} To perform the analyses presented in this section, the model structure \verb|M| has to be created. For more information, see Sec.~\ref{Section:ModelSimNullclines}.

\item \noindent \textit{Read the content of an `.auto' file.} The procedure to read an `.auto' file contains both the $I_0$- and the ($I_0$, $g_K$)-BD is the same as the one presented in Sec.~\ref{Section:1PBD&Eig}. This time, the `Summary' section in the command window presents two diagrams.
\begin{verbatim}
    Summary:
        1P-BD - Name: BD1_i0 - Main: i0
        2P-BD - Name: BD2_i0_gk - Main: i0 - Secondary: gk
\end{verbatim}
The \verb|AR| structure will have the additional BD structure \verb|BD2_i0_gk|.  

\item \noindent \textit{Visualize a 2P-BD.} The next step consists in visualizing the ($I_0$, $g_K$)-BD. In our implementation of system~\ref{SYS:HH}, $I_0$ and $g_K$ are defined as \verb|i0| and \verb|gk|. The function to visualize a 2P-BD is the same used to visualize the 1P-BD: \hyperlink{FUNCTION:Func_VisualizeDiagram}{\texttt{Func\_VisualizeDiagram()}}. This function detects whether the specified BD is a one- or a two-parameter continuation and sets up the visualization automatically. The result is shown in Fig.~\ref{Fig:DEMO3:DSWEB}(A). By default, along the x and y-axes, the main and secondary bifurcation parameters are visualized.
\begin{verbatim}
    fig = figure();
    Func_VisualizeDiagram(M,AR.BD2_i0_gk)
\end{verbatim}
Fig.~\ref{Fig:DEMO3:DSWEB}(B) shows the integration between the $I_0$- and the ($I_0$, $g_K$)-BDs. For this task, we change the axis variables of the $I_0$-BD to \verb|{i0, gk}| by specifying them via the \verb|VAR| field. Alternatively, visualize the BD in the ($I_0$, $g_K$, $V$)-space.

\begin{figure}[ht!]
    \centering
    \includegraphics{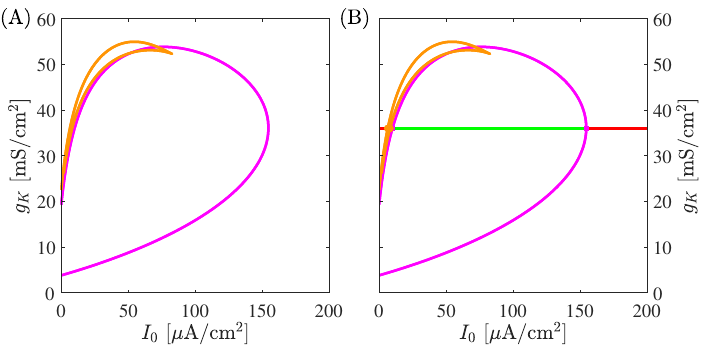}
    \caption{In both panels of the figure, the magenta/orange curves represent the HBs/SNPOs in the ($I_0$, $g_K$) parameter space. (A): ($I_0$, $g_K$)-BD of system (\ref{SYS:HH}). (B) Integration between ($I_0$, $g_K$)-BD of system (\ref{SYS:HH}) and $I_0$-BD. In this panel, the green/red line represents the stable LC/EQ presented in the $I_0$-BD (see Fig.~\ref{Fig:Demo2:1PBD&Eig}).} \label{Fig:DEMO3:DSWEB}
    \vspace*{12pt}
\end{figure}

\end{enumerate}

% Multiple BDs -------------------------------------------------

\subsubsection{Multiple BDs}\label{Section:Example:DEMO3}

\noindent\textbf{Demo}. The accompanying MATLAB file can be found in the folder \verb|DEMOs/DEMO4_HH_MULTI|.

\noindent \textbf{Model.} This section uses the HH model introduced in Sec.~\ref{SYS:HH}.

\noindent \textbf{Data.} The data needed for this task is
\begin{enumerate}[label=(\roman*)]
    \item `.ode' file: File containing an implementation of the HH model~\ref{SYS:HH}.
    \item `.auto' file: File containing the $I_0$-BD for \textit{two} different values of the parameter $g_K$, along with the ($I_0$, $g_K$)-BD. To save these \textit{three} BDs to a single `.auto' file, compute the two $I_0$-BDs, and then create the ($I_0$, $g_K$)-BD by continuing the detected bifurcations along $g_K$. Finally, export the results as specified in Table~\ref{Table:ExportXPPAUT}.
\end{enumerate}

\noindent \textbf{Steps.}
\begin{enumerate}[label=(\arabic*)]
%    \item Visualize a 1P- and a 2P-BD in the same 3D plot
%    \item Visualize multiple 1P-BDs with a 2P-BD in the same 3D plot
%\end{enumerate}

\item \textit{Preliminary step.} To perform the analyses presented in this section, the model structure \verb|M| has to be created. For more information, see Sec.~\ref{Section:ModelSimNullclines}. Moreover, the $I_0$- and the $(I_0,g_K)$-BD has to be loaded using the function \verb|Func_ReadAutoRepo()| (see Sec.~\ref{Section:1PBD&Eig}). 

\item \textit{Visualize a 1P- and a 2P-BD in the same 3D plot.} The function \verb|Func_VisualizeDiagram()| accepts as input the optional field \verb|VAR|. If used for a three-dimensional diagram, this field must be specified as a cell vector with three elements. Here, we want to integrate the information of the two BDs in Figs.~\ref{Fig:Demo2:1PBD&Eig}(A) and~\ref{Fig:DEMO3:DSWEB}(A). \verb|VAR| is set to \verb|{i0,gk,v}|, so that the components visualized along the x/y/z-axes are $I_0$/$g_K$/$V$, respectively. The result for the following code is presented in Fig.~\ref{Fig:DEMO4}.
\begin{verbatim}
    fig = figure();
    Func_VisualizeDiagram(M,AR.BD2_i0_gk,`VAR',{`i0',`gk',`v'})
    Func_VisualizeDiagram(M,AR.BD1_i0   ,`VAR',{`i0',`gk',`v'})
    Func_VisualizeLabPoints(M,AR.BD1_i0 ,`VAR',{`i0',`gk',`v'})
\end{verbatim}

\noindent Note that all the visualization functions use the \verb|BD| structure. An alternative to the current visualization is presented in the next section.

\item \textit{Visualize multiple 1P-BDs and a 2P-BD in the same 3D plot.} This can be achieved by adding a 1P-BD (from the same or another \verb|AR| structure) to a figure with a 1P- and 2P-BD, using the functions \verb|Func_VisualizeDiagram()| and \verb|Func_VisualizeLabPoints()|. We omit the results from this section here but encourage the reader to explore the results interactively in the accompanying `.m' demo files.

\begin{figure}[ht!]
    \centering
    \includegraphics{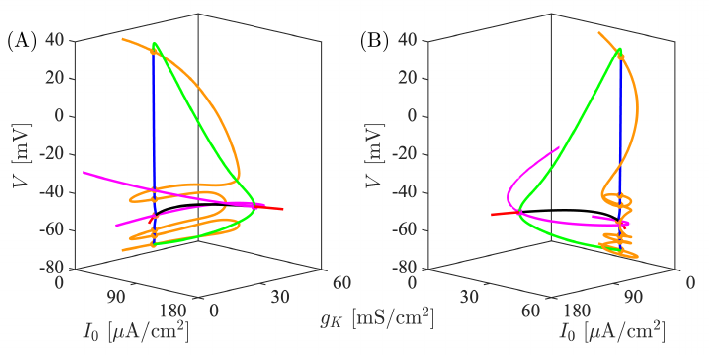}
    \caption{Union of 1P-BD from Fig.~\ref{Fig:Demo2:1PBD&Eig}(A) and 2P-BD from Fig.~\ref{Fig:DEMO3:DSWEB}(A). Both panels visualize the same content, but the right has reversed x- and y-axes. The \textit{magenta curve} connects the two HBs, hence it represents the continuation of these bifurcations in the two-parameter space. The \textit{orange curve} represents the two-parameter continuation of the SNPOs.}
    \label{Fig:DEMO4}
    \vspace*{12pt}
\end{figure}

\textit{Additional Note: 3D plot.} Fig.~\ref{Fig:DEMO4} illustrates the visualization of a 1P-BD in a 3D plot with respect to a single dynamical variable and two parameters. Instead, one can also visualize a 1P-BD with respect to two dynamical variables and one parameter. When a 1P-BD has stable limit cycles, they are represented by branches of the maxima and minima for each value of the bifurcation parameter. However, a discrepancy can arise between the minimum and maximum branches of the 1P-BD and the reconstruction of the manifold of limit cycles (see Fig.~\ref{Fig:DEMO6} for reconstruction). This occurs because the maxima and minima are only evaluated in the visualized dynamical variable, but the maxima and minima in one variable will not necessarily coincide with the extrema of the other dynamical variables.

\end{enumerate}

%%%%%%%%%%%%%%%%%%%%%%%%%% Averaging %%%%%%%%%%%%%%%%%%%%%%%%%%

\subsection{Averaging}\label{Section:Averaging}

\noindent\textbf{Demo}. The accompanying MATLAB file can be found in the folder \verb|DEMOs/DEMO5_CK_AVG|.

\noindent \textbf{Model.} We consider the system of ODEs presented in~\cite{CK1983}, modeling bursting electrical activity in pancreatic beta cells. The ODEs governing the behavior of the dynamical variables are the following,
\begin{equation}
    \begin{aligned}
        c_m \dot V &= -I_K(V,n) - I_{Ca}(V) - I_{K(\text{Ca})}(V,c) - I_{K(\text{ATP})}(V), \\
        \dot n &= (n_{\infty}(V) - n) \tau_n^{-1}, \\
        \dot c &= -f(\alpha I_{\text{Ca}}(V) + k_{\text{PMCA}}c).
    \end{aligned}
    \label{SYS:CK}
\end{equation}
where $V$ is the cell's membrane potential, $n$ is the activation of K\textsuperscript{+} channel, while $c$ represents the intracellular calcium concentration. The fraction of free calcium in $\beta$-cells is denoted with $f$. $I_K$, $I_{K(Ca)}$, $I_{K(ATP)}$, and $I_{Ca}$ represent the K\textsuperscript{+}, Ca\textsuperscript{2+}-activated K\textsuperscript{+}, ATP-sensitive K\textsuperscript{+} and Ca\textsuperscript{2+} currents, respectively. This system can generate bursting activity. One approach to unveil the dynamical mechanism leading to the generation of this phenomenon is to consider $V$ and $n$ as fast, while $c$ as a slow variable, as initially presented in~\cite{Rinzel1985}. For this reason, we can study the full evolution of the system, treating $c$ as a fixed parameter and investigating the stable and unstable attractors of the fast subsystem in ($V$, $n$). The fast subsystem is defined as
\begin{equation}
    \begin{aligned}
        c_m \dot V &= -I_K(V,n) - I_{Ca}(V) - I_{K(\text{Ca})}(V|c) - I_{K(\text{ATP})}(V), \\
        \dot n &= (n_{\infty}(V) - n) \tau_n^{-1}.
    \end{aligned}
    \label{SYS:CK:FAST}
\end{equation}
To further elucidate the bursting mechanism, the averaging theory is applied~\cite{Bertram1995}. This type of analysis can be performed by using XPPLORE and XPPAUT.

\noindent \textbf{Data.} The data required to perform averaging analyses are
\begin{enumerate}[label=(\roman*)]
    \item `.ode' file: File containing an implementation of system (\ref{SYS:CK:FAST}).
    \item `.auto' file: File containing the $c$-BD of system (\ref{SYS:CK:FAST}). To export the periodic orbits over the branch of stable LCs emanating from the supercritical HB, the parameter \verb|NPR| under the \verb|nUmerics| menu of XPPAUT must be set to 1. This setting controls how often all the information related to a fixed point or a periodic orbit is exported. Through the approximation of the periodic orbits, the averaging theory will be applied numerically.
\end{enumerate}

\noindent \textbf{Steps.}
\begin{enumerate}[label=(\arabic*)]
%    \item Extract the periodic orbits
%    \item Applying averaging
%    \item Visualization of 1P-BD with averaging results
%    \item Convert \& export `.auto' to `.dat' file
%\end{enumerate}
\item \textit{Preliminary step.} The analyses requires the creation of the model structure \verb|M| (see Sec.~\ref{Section:ModelSimNullclines}) and the loading of the $c$-BD (see Sec.~\ref{Section:1PBD&Eig}). 

\item \textit{Extract the periodic orbits.} To perform averaging analysis, we need to extract the LC from the \verb|AR| structure. We can retrieve all the special trajectories from the $c$-BD by using the function \verb|Func_GetTRJ()|.
\begin{verbatim}
    TRJ = Func_GetTRJ(M,AR.BD1_c);
\end{verbatim}
The variable \verb|TRJ| is a structure of special trajectories. For more information related to this structure, please visit Sec.~\ref{Section:Trajectories}. By typing `\verb|TRJ|' in the MATLAB command window, the following output will be presented.
\begin{verbatim}
    TRJ = 
        struct with fields:
          TRJ1: [1×1 struct]
          TRJ2: [1×1 struct]
          TRJ3: [1×1 struct]
          ...
        TRJ366: [1×1 struct]
          nTRJ: 366
\end{verbatim}
\verb|TRJ| is a structure whose fields are the special trajectories \verb|TRJi| with \verb|i| in \verb|{1,...,nTRJ}|. \verb|nTRJ| is the number of labeled points associated with a periodic orbit stored under the field \verb|LABPTs| of \verb|BD1_c|. 
%Please refer to the accompanying \textit{.m} file to further explore the \verb|TRJ| structure.

\item \textit{Applying averaging.} %We are ready to apply averaging.
Averaging theory consists of calculating the average variation of the slow variable over the numerically approximated LC. In the following, we will refer to the LC with $\gamma$. Mathematically,
\begin{equation}
    \langle\dot{c}\rangle = -\frac{1}{T_{\gamma}}\int_{0}^{T_\gamma}{[\alpha I_{Ca}(V_{\gamma}(t)) + k_{PMCA}c]}\:dt.
    \label{SYS:CK:AVG}
\end{equation}
$V_\gamma(t)$ is the $V$ coordinate of the fast subsystem during its evolution along the stable LC $\gamma$. $T_\gamma$ is the period of the LC, $c$ is the parameter value for which $\gamma$ exists and $\langle \dot{c} \rangle$ is the average variation of the slow variable over $\gamma$.
To compute the previous integral, we can create a custom function which cycles through all the special trajectories, taking advantage of the organization of the \verb|TRJ| structure, and computes the previous integral.
\begin{verbatim}
    function [c,J,BZ] = Func_Averaging(M,TRJ,kpmca)

    % FUNCTIONs
    mI   = @(V)       0.5*(1+tanh((V-M.P.vm)/M.P.sm));
    ICa  = @(V)       M.P.gca*mI(V).*(V-M.P.vca);
    Jmem = @(c,V,K) -(M.P.alpha*ICa(V) + K*c);
    
    % INITIALIZATION
    J = zeros(TRJ.nTRJ,1);
    c = zeros(TRJ.nTRJ,1);
    F = fieldnames(TRJ);
    
    % AVERAGING
    for iTRJ = 1:1:TRJ.nTRJ
        t = TRJ.(F{iTRJ}).PO.t;
        V = TRJ.(F{iTRJ}).PO.v;
    
        T       = TRJ.(F{iTRJ}).P.T;
        c(iTRJ) = TRJ.(F{iTRJ}).P.c;
        J(iTRJ) = (1/T)*trapz(t*T,Jmem(c(iTRJ),V,kpmca));
    end
    
    % EXTRACTION
    [~,IZ] = min(abs(J));
    BZ = TRJ.(F{IZ}).PTi;
    
    end
\end{verbatim}
In the `functions' section of the code, the functions needed to perform averaging as specified in Eq. (\ref{SYS:CK:AVG}) are defined. In the `initialization' section, the vectors and the variables used to keep track of the computed values are initialized. Finally, in the `averaging' section, the code cycles through all of the trajectories, and the integral in Eq. (\ref{SYS:CK:AVG}) is numerically approximated. Lastly, in the `extraction' section, the name of the labeled point associated with the orbit with the lowest absolute average variation is saved in the variable \verb|BZ|. The output values of this custom function are \verb|c|, \verb|J|, and \verb|BZ|. \verb|c| and \verb|J| are two one-dimensional vectors of the same length containing the value of parameter $c$ associated with the LC and the average variation $\langle \dot{c} \rangle$ computed along them.

\item\textit{Visualization of 1P-BD with averaging results.} The results of the averaging theory can be presented in different ways. In this file, we present the slow average variation ($\langle \dot{c}\rangle$) in correspondence with different periodic orbits alongside the $c$-BD of the fast subsystem. The results are shown in Fig.~\ref{Fig:DEMO5}. To create this image, we use the functions \verb|Func_VisualizeDiagram()| and \verb|Func_VisualizeLabPoints()| to handle the visualization of the $c$-BD. The remaining parts of the image are handled through simple default MATLAB commands. A compressed version of the code is presented in the following.
\begin{verbatim}
    fig = figure();
    tiledlayout(1,2,`TileSpacing',`Compact',`Padding',`Compact')

    nexttile()
    Func_VisualizeDiagram(M,AR.BD1_c)
    Func_VisualizeLabPoints(M,AR.BD1_c)

    hold on
    plot(AR.BD1_c.LABPTs.(BZ).c,AR.BD1_c.LABPTs.(BZ).vU,...
    `Color',`m',`Marker',`o',`MarkerFaceColor',`m',`MarkerSize',3)
    plot(AR.BD1_c.LABPTs.(BZ).c,AR.BD1_c.LABPTs.(BZ).vL,...
    `Color',`m',`Marker',`o',`MarkerFaceColor',`m',`MarkerSize',3)
    hold off

    nexttile()
    plot(c,J,`Color',`b',`LineWidth',1.2)

    hold on
    plot(c(islocalmin(abs(J))),J(islocalmin(abs(J))),...
    `Color',`m',`Marker',`o',`MarkerFaceColor',`m',`MarkerSize',3)
    hold off
\end{verbatim}

\begin{figure}[ht!]
    \centering
    \includegraphics{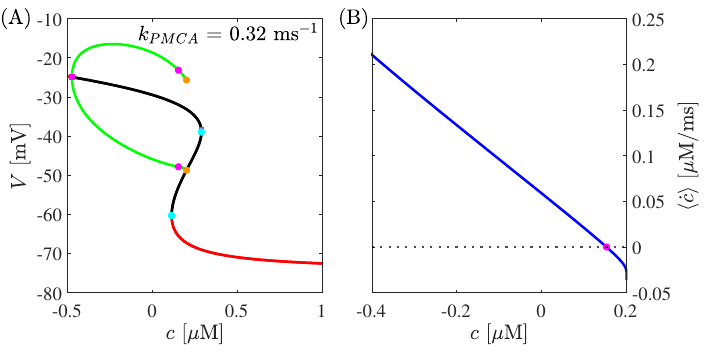}
    \caption{(A) The $c$-BD of the fast subsystem and the results of the averaging theory. The black and red curves represent unstable and stable fixed points, respectively. The green curves are the maxima and minima of stable limit cycles. The magenta squares are HBs, and the cyan and orange dots are SN and SNPO, respectively. The left-most HB is supercritical; in fact, the limit cycles emanating from here are stable (green curves). The magenta dots over the branch of stable limit cycles represent the minima and maxima of the PO, presenting the numerically approximated $\langle \dot{c} \rangle$ closest to 0. (B) The dependency of $\langle \dot{c} \rangle$ on $c$ is visualized. The magenta dot corresponds to the one in panel (A), hence it is the solution with an average variation closest to zero. For this analysis, $k_{PMCA}$ was set to 0.32 ms\textsuperscript{-1}.}
    \label{Fig:DEMO5}
    \vspace*{12pt}
\end{figure}

\item \textit{Convert and export the `.auto' in a `.dat' file.} To simplify the conversion of an `.auto' file into an XPPAUT-interpretable file, we use the function \verb|Func_WritePoints()|, which requires as input the model structure \verb|M|, the bifurcation diagram to be exported and the new BD file's name. 
\begin{verbatim}
    Func_WritePoints(M,AR.BD1_i0,`BD.dat');
\end{verbatim}
The function exports results to `BD.dat', which can be loaded in XPPAUT through \verb|Graphic stuff| $\rightarrow$ \verb|Freeze| $\rightarrow$ \verb|Bif.Diag| in XPPAUT.
\end{enumerate}

%%%%%%%%%%%%%%%%%%%%%%%%%% Slow Manifold %%%%%%%%%%%%%%%%%%%%%%%%%%

\subsection{Reconstructing manifolds}

\subsubsection{Manifold of limit cycles}

\noindent\textbf{Demo}. The accompanying MATLAB file can be found in the folder \verb|DEMOs/DEMO6_CK_MLC|.

\noindent \textbf{Model.} This demo considers the dynamical system in Eq. (\ref{SYS:CK}).

\noindent \textbf{Data.} The data required to reconstruct the stable limit cycles of system (\ref{SYS:CK}) are the following,
\begin{enumerate}[label=(\roman*)]
    \item `.ode' file: File containing the implementation of system (\ref{SYS:CK}).
    \item `.auto' file: File with the results of the numerical continuations ($c$-BD).
\end{enumerate}

\noindent \textbf{Steps.}

\begin{enumerate}[label=(\arabic*)]
%    \item Reconstruct the manifold
%    \item Visualize the manifold
%\end{enumerate}

\item \textit{Preliminary step.} See the `Preliminary Step' in Sec.~\ref{Section:Averaging}. Moreover, we assume that the user has loaded the `.auto' file and extracted the calculated special trajectories TRJ, as explained in Sec.~\ref{Section:Averaging}.

\item \textit{Reconstruct the surface of limit cycles.}  To convert the continuation output into a smooth surface in MATLAB, we create a custom function \verb|Func_Manifold()|. This leverages the uniform length of all trajectories computed via numerical continuation, as they possess an equal number of time points. The rationale of the function presented below is simple. The trajectories are treated as column vectors, which are concatenated horizontally. This procedure allows us to create matrices containing a numerical, approximated reconstruction of the surface of limit cycles within the bifurcation diagram.
\begin{verbatim}
    function S = Func_Manifold(TRJ,VAR,[PAR])

    % INITIALIZATION
    if nargin == 2     , nPAR = 0;
    elseif isempty(PAR), nPAR = 0;
    else               , nPAR = length(PAR);
    end
    
    nS   = TRJ.nTRJ;
    nTS  = length(TRJ.TRJ1.PO.t);
    nVAR = length(VAR);
    for iVAR = 1:nVAR, S.(VAR{iVAR}) = zeros(nTS,nS); end
    for iPAR = 1:nPAR, S.(PAR{iPAR}) = zeros(nTS,nS); end
    
    % CONCATENATION
    for iT = 1:nS
        TR = sprintf(`TRJ%i',iT); 
        PO = TRJ.(TR).PO; 
        P  = TRJ.(TR).P; 
        for iVAR = 1:1:nVAR, S.(VAR{iVAR})(:,iT) = PO.(VAR{iVAR}); end
        for iPAR = 1:1:nPAR, S.(PAR{iPAR})(:,iT) = ones(nTS,1).*P.(PAR{iPAR}); end
    end
\end{verbatim}
The processing of the data is based on two steps. First, the matrices containing the limit cycles are `initialized'. During the second step, the trajectories are `concatenated'. Note that the dynamical variables of the fast subsystem are treated separately from the bifurcation parameter $c$ here.

\item \textit{Visualize the surface.} To visualize the reconstructed stable limit cycle manifold and the bifurcation diagram, we use the built-in MATLAB function \verb|surf()| and the XPPLORE function \verb|Func_VisualizeDiagram()|. For the visualization of a sample trajectory in the three-dimensional BD, we use \verb|plot()|. The results are shown in Fig.~\ref{Fig:DEMO6}.\\
See Sec.~\ref{Section:Example:DEMO3} for a note regarding the reconstructed manifold of limit cycles and the calculation of branches representing limit cycles in 1P-BDs.

\begin{figure}[ht!]
    \centering
    \includegraphics{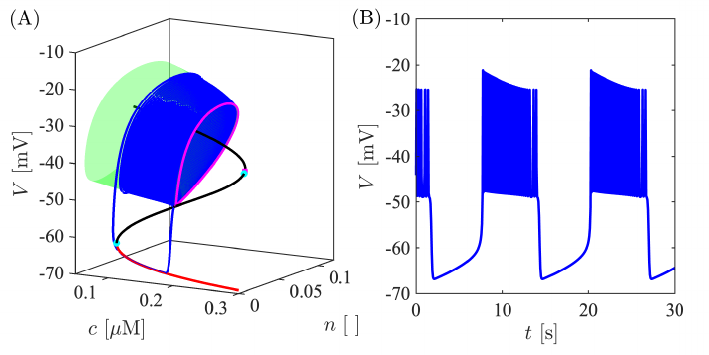}
    \caption{(A) Numerical reconstruction of the manifold of stable limit cycles of system (\ref{SYS:CK}). The magenta curve is a numerical approximation of the homoclinic trajectory. The blue curve is a trajectory of system (\ref{SYS:CK}). The green surface is the stable manifold of limit cycle. The s-shape red/black curve represents the stable/unstable fixed point of the system in Eq. (\ref{SYS:CK:FAST}). It corresponds to the curve of stable/unstable points presented in Fig.~\ref{Fig:DEMO5}. The cyan circles are SN bifurcations. (B) Temporal evolution of the trajectory presented in panel (A).}
    \label{Fig:DEMO6}
\end{figure}

\end{enumerate}

\subsubsection{Slow manifold}

\noindent\textbf{Demo}. The accompanying MATLAB file can be found in the folder \verb|DEMOs/DEMO7_FHN_SM|.

\noindent \textbf{Model.} The system of ODEs considered for this demo is the silent phase of a self-coupled FitzHugh-Nagumo model presented in~\cite{DKO2008CHAOS}.
\begin{equation}
    \begin{aligned}
        \dot{v} &= h-0.5(v^3-v+1)-\gamma s v \\
        \dot{h} &= -\epsilon(2h + 2.6v) \\
        \dot{s} &= -\epsilon(\delta s)
    \end{aligned}
    \label{SYS:FHN:3D}
\end{equation}

\noindent \textbf{Data.} The data required to reconstruct the attracting and the repelling slow manifold of system (\ref{SYS:FHN:3D}) are the following,
\begin{enumerate}[label=(\roman*)]
    \item `.ode' file: File containing the implementation of system (\ref{SYS:FHN:3D}) to reconstruct the repelling slow manifold with a two-point BVP approach.
    \item `.auto' file: File with the results of the numerical continuations.
\end{enumerate}

\noindent \textbf{Steps.}
\begin{enumerate}[label=(\arabic*)]
%    \item Reconstruct the slow manifold
%    \item Visualize the slow manifold
%\end{enumerate}

\item \textbf{Additional notes.} The procedure illustrated here does not provide any detail regarding the method to reconstruct the slow manifold with the usage of XPPAUT. The general idea exploited here is based on solving a two-point BVP using a homotopy method with two steps. The first step in the homotopy continues a trivial solution (the folded node) away from the folded node plane ($\Sigma$); instead, the second step aims to gain normal hyperbolicity by continuing the solution along the attracting or repelling submanifolds of the critical manifold. For technical information on how this procedure works in XPPAUT, we ask the reader to refer to future publications. For more details regarding the formalization of this approach, we recommend the reader the pioneering work presented in \cite{DKO2008CHAOS, DKO2008SIAM, DKO2010}. Since no details regarding the calculation of the slow manifolds are provided in this paper, a downsampled version of the dataset containing the results of the continuation can be found in the folder associated with this example. 

\item \textit{Preliminary step.} This demo requires the creation of the model structure \verb|M| (see Sec.~\ref{Section:ModelSimNullclines}), the loading of the continuation results stored in an `.auto' file in an \verb|AR| structure (see Sec.~\ref{Section:1PBD&Eig}) and the extraction of all the approximated solutions into the \verb|TRJ| structure (see Sec.~\ref{Section:Averaging}).

To perform the analyses presented in this section, the model structure \verb|M| has to be created. For more information, see Sec.~\ref{Section:ModelSimNullclines}. Moreover, the results of the numerical continuation have to be loaded using the function \verb|Func_ReadAutoRepo()| (see Sec.~\ref{Section:1PBD&Eig}). 

\item \textit{Reconstruct the slow manifold.} Similarly to the previous manifold reconstruction, we assume that the user has loaded the `.auto' file and extracted the calculated special trajectories TRJ. We use the custom function \verb|Func_Manifold()| introduced for limit cycles to reconstruct the slow manifold. The projection of the slow manifold over the terminal plane are extracted through the function \verb|Func_SlowManifoldProjection()|. The rationale of this function is the same as \verb|Func_Manifold()|, with the only difference being that, instead of concatenating all the trajectories, we only concatenate their last points.

\item \textit{Visualize the slow manifold.} To visualize the reconstructed slow manifolds, we use the built-in MATLAB function \verb|surf()|. For the visualization of the projection of the attracting and the repelling slow manifold over the terminal plane, we use the command \verb|plot()|. The results are shown in Fig.~\ref{Fig:DEMO7}.

\begin{figure}[ht!]
    \centering
    \includegraphics{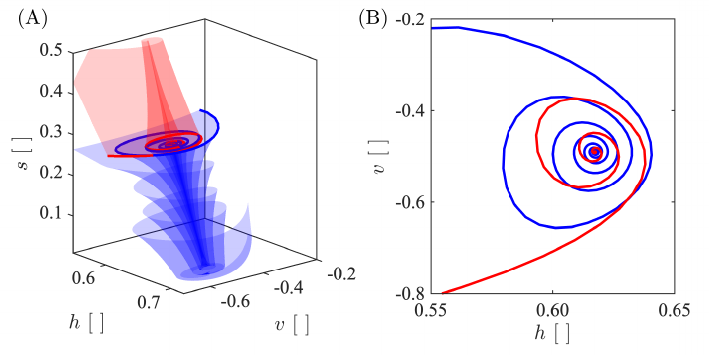}
    \caption{(A) The red/blue surface is a numerical reconstruction of the attracting/repelling slow manifold of system (\ref{SYS:FHN:3D}). (B) The red/blue curve represents the projection of the attracting and repelling slow manifold over the terminal plane $\Sigma$. This figure is a partial downsampled version of Fig. 3 presented in \cite{DKO2008CHAOS}.}
    \label{Fig:DEMO7}
    \vspace*{12pt}
\end{figure}

\end{enumerate}

\newpage
\section{Discussion and conclusions}
After a review of the scientific literature and the evaluation of the software packages to investigate ODE-based dynamical systems, we developed XPPLORE to increase the accessibility of continuation results for XPPAUT users. XPPLORE preserves the usual XPPAUT workflows and leverages the continuation calculations exported in `.auto' files. The user-friendly functions in our toolbox allow the straightforward interaction between the user and the XPPLORE data structures to manipulate, explore, and visualize the continuation results. This software opens the door to multiple types of analyses, requiring results produced through numerical continuations, extending the XPPAUT capabilities. We first showed that XPPLORE can perform simple tasks, performed by existing extension packages, such as the visualization of continuation results, e.g. 1P/2P-BD, nullclines, and simulations. In addition, we presented how our package can be used (1) to directly manipulate continuation calculations to obtain more advanced results such as the integration of one- and two- parameter bifurcation diagrams in two- and three-dimensional space, (2) to apply averaging theory, (3) to reconstruct different types of manifolds. %manifold of stable limit cycles in~\cite{CK1983}, and finally (4) to reconstruct the slow manifold close to a node in a silent formulation of the FitzHugh-Nagumo model presented in~\cite{CK1983}. 

We are confident that the functionalities of our packages will allow the scientific community to improve the use of XPPAUT. We believe that our toolbox can be used for analyses well beyond the examples presented in this paper, e.g. it can be used to initialize BVPs to reconstruct super-slow nullclines in systems with two or more super-slow variables. Moreover, another use can be the extraction of the calculations of isochrones and isostable curves in two-dimensional dynamical systems. Due to the user-friendly design and the flexibility through its implementation in MATLAB, we believe that this package can easily be translated into other programming languages. In conclusion, XPPLORE sheds new light on hidden XPPAUT capabilities, paving the way to new analyses and exploration.

\section{Acknowledgements}

The authors thank Jonathan Rubin and Morten Gram Pedersen for their valuable suggestions and discussions related to this work. Matteo Martin gratefully acknowledges the support of the Fondazione Aldo Gini of Padova for this work. Anna Kishida Thomas received support from NIH R01NS125814.
%I will need to ask Jon if I should add my NIH grant also.%

\newpage
\section{Appendix}%\label{Section:Appendix}

\subsection*{Function documentation}

\noindent In the following, we present the functions implemented in this package, their mandatory and optional inputs, as well as their outputs. Table~\ref{Table_Func:Functions} lists the categories of functions available, the function names, and a hyperlink to their detailed documentation. The functions are documented in alphabetical order.

\begin{table}[ht]
    \centering
    \caption{This table contains all of the available functions in this package, sorted by category, with references to their full documentation. \label{Table_Func:Functions}}
    {\tabcolsep13pt\begin{tabular}{cc}\\[-2pt]
        \toprule
        Category & Function \\[6pt]
        \hline
        \multirow{4}{*}{Reading} & \hyperlink{FUNCTION:Func_ReadModel}{\texttt{Func\_ReadModel()}}\\
         & \hyperlink{FUNCTION:Func_ReadData()}{\texttt{Func\_ReadData()}}\\
         & \hyperlink{FUNCTION:Func_ReadNullclines}{\texttt{Func\_ReadNullclines()}} \\
         & \hyperlink{FUNCTION:Func_ReadAutoRepo}{\texttt{Func\_ReadAutoRepo()}} \\
         \multirow{4}{*}{Visualization} & \hyperlink{FUNCTION:Func_VisualizeDiagram}{\texttt{Func\_VisualizeDiagram()}}\\
         & \hyperlink{FUNCTION:Func_VisualizeLabPoints}{\texttt{Func\_VisualizeLabPoints()}}\\
         & \hyperlink{FUNCTION:Func_VisualizeEig}{\texttt{Func\_VisualizeEig()}}\\
         & \hyperlink{FUNCTION:Func_VisualizeNullclines}{\texttt{Func\_VisualizeNullclines()}}\\
         
         \multirow{3}{*}{Extraction} & \hyperlink{FUNCTION:Func_GetEig}{\texttt{Func\_GetEig()}}\\
         & \hyperlink{FUNCTION:Func_GetTRJ}{\texttt{Func\_GetTRJ()}} \\
         & \hyperlink{FUNCTION:Func_WritePoints}{\texttt{Func\_WritePoints()}} \\
         
         \multirow{2}{*}{Options} & \hyperlink{FUNCTION:Func_DOF}{Func\_DOF()}\\
         & \hyperlink{FUNCTION:Func_DOBD}{\texttt{Func\_DOBD()}} \\
         
         \multirow{2}{*}{Figures} & \hyperlink{FUNCTION:Func_FigExport}{\texttt{Func\_FigExport()}} \\
         & \hyperlink{FUNCTION:Func_FigStyle}{Func\_FigStyle()}\\
         \hline
    \end{tabular}}
\end{table}

\vspace{0.7cm}

\begin{sidewaysfigure}
    \centering
    \centerline{\includegraphics[width=9in]{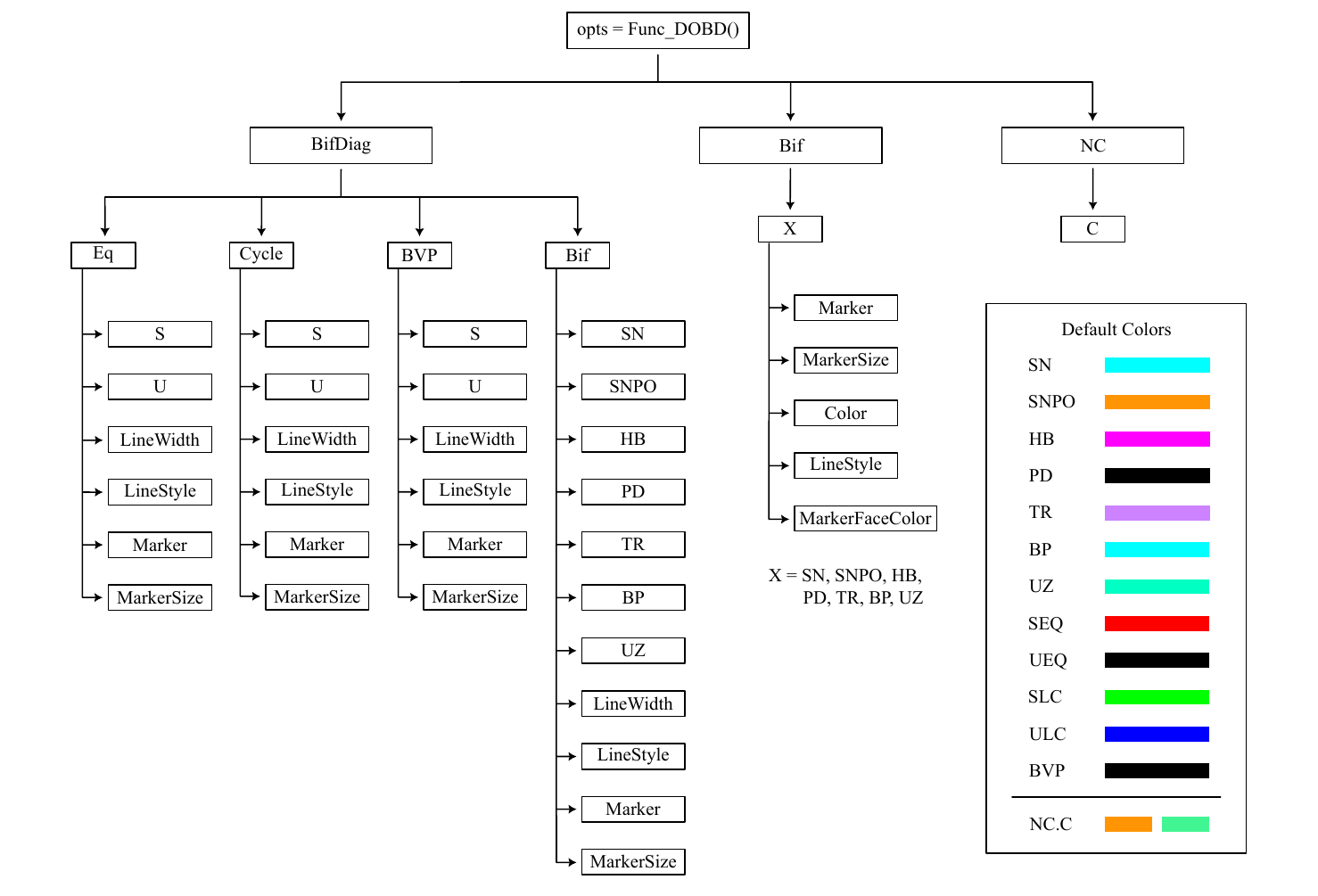}}
    \caption{Graphical representation of options structure \texttt{opts} for bifurcation diagrams. The inset shows default color choices for bifurcation points and branches within bifurcation diagrams. The options structure consists of three sub-structures: \texttt{BifDiag}, containing settings for branches within 1P- and 2P-BDs; \texttt{Bif}, containing specifications regarding bifurcation point appearances; \texttt{NC}, the figure settings for nullclines. Note that \texttt{opts.BifDiag.Bif} target the appearance of branches of bifurcations in a 2P-BD while \texttt{opts.Bif} refers to bifurcation points within 1P-BDs. \textit{S=stable, U=unstable, BP=bifurcation point, C=color.}}
    \label{FIG:OPTSBD}
\end{sidewaysfigure}

%%%% FUNC_DOBD()
\noindent \hypertarget{FUNCTION:Func_DOBD}{\texttt{opts = Func\_DOBD()}}%\label{Func_DOBD}
    \begin{enumerate}[label=(\roman*)]
    \item \verb|DESCRIPTION|: \\
    This function creates the default structure containing the visualization settings for the bifurcation diagrams, labeled points, and nullclines. All the inputs to this function are optional and must be passed as name-value pairs. The accepted inputs are RGB vectors specified in the MATLAB format. Hence, vectors of three numbers between 0 and 1.
    \item \verb|INPUTs|: 
    \begin{enumerate}[label=(\alph*)]
    \item \verb|SNPO|: Color used for SNPO.
    \item \verb|SN|: Color used for SN.
    \item \verb|TR|: Color used for TR.
    \item \verb|HB|: Color used for HB.
    \item \verb|PD|: Color used for PD.
    \item \verb|BP|: Color used for BP.
    \item \verb|UZ|: Color used for UZ.
    \item \verb|SEQ|: Color used for stable EQ point.
    \item \verb|UEQ|: Color used for unstable EQ point.
    \item \verb|SLC|: Color used for stable LC.
    \item \verb|ULC|: Color used for unstable LC.
    \item \verb|BVP|: Color used for BVP.
    \end{enumerate}
    \item \verb|OUTPUTs|: \\
    The output is a structure of different types of visualization settings used for the visualization of bifurcation diagrams, labeled points, and nullclines.
    \item \verb|ADDITIONAL NOTEs|: \\
    We encourage the user to edit this file to customize the visualization. Figure~\ref{FIG:OPTSBD} contains a complete list of these settings.
    \end{enumerate}

\vspace{0.7cm}

%%%% FUNC_DOF()
\noindent \hypertarget{FUNCTION:Func_DOF}{\texttt{opts = Func\_DOF()}}
    \begin{enumerate}[label=(\roman*)]
    \item \verb|DESCRIPTION|: \\
    This function creates a structure containing the settings for the style of the figure. All of the inputs to this function are optional and must be passed as name-value pairs.
    \item \verb|INPUTs|:
    \begin{enumerate}[label=(\alph*)]
    \item \verb|format|: This field accepts a sequence of characters. The valid formats are specified in the \verb|formattype| field of the \verb|print()| function implemented in MATLAB. The values determine the format of the exported image.
    \item \verb|extension|: This field accepts a sequence of characters representing the extension of the exported file. Keep in mind that, if the user changes the \verb|format| field of this structure, they should change this field as well.
    \item \verb|resolution|: In this field, the user may specify the resolution in dpi of the image through a sequence of characters. Specifically, the entry must have the format \verb|`-rXXX'|, where \verb|`XXX'| is the resolution in dpi. For more information regarding this field, please visit the official MATLAB documentation of the function \verb|print()|. 
    \item \verb|units|: Figure's unit of measure. This field changes the unit of measure in which all the options specified in this structure will be interpreted when the style is applied to the figure.
    \item \verb|width|: Width of the exported figure in \verb|units|. 
    \item \verb|height|: Height of the exported figure in \verb|units|.
    \item \verb|Layer|: It defines the placement of grid lines and tick marks compared to the visualized graphic objects.
    \item \verb|ClippingStyle|: Character array representing the type of clipping applied to the figure. 
    \item \verb|FontName|: Name of the axes and title font.
    \item \verb|FontSize|: Size of the text visualized in the current axes.
    \item \verb|gridLineStyle|: Style of the grid line.
    \item \verb|minorGridLineStyle|: Style of the minor grid line. 
    \item \verb|gridAlpha|: Alpha levels of the grid line. 
    \item \verb|minorGridAlpha|: Alpha levels of the minor grid line.
    \end{enumerate}
    \item \verb|OUTPUTs|: \\
    The only output of this function is the option structure \verb|opts| containing all the different fields specified previously.
    \item \verb|ADDITIONAL NOTEs|:
    \begin{enumerate}[nolistsep]
    \item We strongly encourage the user to edit the default values implemented in this function. 
    \item To understand how the \verb|ClippingStyle| and \verb|Layer| fields affect the visualization, the user should refer to the MATLAB documentation.
    \item When the user changes the \verb|units| field of this function, all declared properties are not converted into the new unit of measure automatically. The user must change these manually.
    \item For more information on \verb|gridLineStyle|, \verb|minorGridLineStyle|, \verb|gridAlpha|, and \\
    \verb|minorGridAlpha|, please refer to the MATLAB documentation of the axes objects.
    \item The default values of the previous fields are presented in table~\ref{Table_OPTS:Figure}.
    \begin{table}[ht!]
        \centering
        \caption{This table shows the default values stored in the options structure. \label{Table_OPTS:Figure}}
        {\tabcolsep13pt\begin{tabular}{cccc}\\[-2pt]
            \toprule
            Properties & Default values & Properties & Default values \\
            \hline
            format & `-dpdf' & extension & `.pdf' \\
            resolution & `-r600' & units & `centimeters' \\
            width & 6 & height & 6 \\
            ClippingStyle & rectangle & Layer & top \\
            FontName & Times New Roman & FontSize & 9 \\
            gridLineStyle & `none' & minorGridLineStyle & `none'\\
            gridAlpha & 0 & minorGridAlpha & 0 \\
            \hline
        \end{tabular}}
    \end{table}
    \end{enumerate}
    \end{enumerate}

\vspace{0.7cm}

%%%% FUNC_FIGEXPORT()
\noindent \hypertarget{FUNCTION:Func_FigExport}{\texttt{Func\_FigExport(fig, fileName, "OPTIONs")}}
    \begin{enumerate}[label=(\roman*)]
    \item \verb|DESCRIPTION|: \\
    This function is used to export the figure created through MATLAB in the desired format specified inside the option structure created using \hyperlink{FUNCTION:Func_DOF}{\texttt{Func\_DOF()}}.
    \item \verb|INPUTs|:
    \begin{enumerate}[label=(\alph*)]
    \item \verb|fig|: Object referring to the figure to be exported.
    \item \verb|fileName|: File's name where the exported figure will be saved.
    \item \verb|OPTIONs|: This is an optional field. It has to be specified through a name-value declaration. As input, it accepts an option structure created through the function \hyperlink{FUNCTION:Func_DOF}{\texttt{Func\_DOF()}}.
    \end{enumerate}
    \item \verb|OUTPUTs|: -
    \end{enumerate}

\vspace{0.7cm}

%%%% FUNC_FIGSTYLE()
\noindent \hypertarget{FUNCTION:Func_FigStyle}{\texttt{Func\_FigStyle(fig, "OPTIONs")}}
    \begin{enumerate}[label=(\roman*)]
    \item \verb|DESCRIPTION|: \\
    This function allows the user to apply the style specified in the option structure created through the function \hyperlink{FUNCTION:Func_DOF}{\texttt{Func\_DOF()}} to the input figure. If the only input is a figure object, then the applied style is the default one created through the aforementioned function.
    \item \verb|INPUTs|:
    \begin{enumerate}[label=(\alph*)]
    \item \verb|fig|: This input is the figure object to which the user wants to apply the figure's style.
    \item \verb|OPTIONs|: This is an optional field. It must follow a name-value declaration. It accepts as input an option structure created through the function \hyperlink{FUNCTION:Func_DOF}{\texttt{Func\_DOF()}}.
    \end{enumerate}
    \item \verb|OUTPUTs|: -
    \end{enumerate}  

\vspace{0.7cm}

%%%% FUNC_GETEIG()
\noindent \hypertarget{FUNCTION:Func_GetEig}{\texttt{EIG = Func\_GetEig(BD,"SORTED")}}
    \begin{enumerate}[label=(\roman*)]
    \item \verb|DESCRIPTION|:\\
    This function allows the user to retrieve the eigenvalues along the branches within a bifurcation diagram in a matrix format, for ease of use. Note that these values represent eigenvalues in the case of fixed points and Floquet multipliers in the case of limit cycles. More importantly, for fixed points, the eigenvalues are post-processed with an exponential transformation in XPPAUT, directly. This function is also utilized by \hyperlink{FUNCTION:Func_VisualizeEig}{\texttt{Func\_VisualizeEig()}}.
    \item \verb|INPUTs|:
    \begin{enumerate}[label=(\alph*)]
    \item \verb|BD|: Bifurcation diagram structure.
    \item \verb|"SORTED"|: This is an optional argument which is set to the boolean value \verb|true| by default and indicates whether the eigenvalues of each point within the BD are sorted in ascending order with respect to their real components. Imaginary parts of eigenvalues are ordered accordingly. This argument must be passed as a name-value argument.
    \end{enumerate}
    \item \verb|OUTPUTs|: \\
    The output of this function is a matrix where the first index indicates the branch number (same numbering as in the input BD struct), the second index identifies the point number within that branch, and the third and fourth indices refer to real and imaginary components of the eigenvalue, respectively.
    \end{enumerate}   

\vspace{0.7cm}

%%%% FUNC_GETTRJ()
\noindent \hypertarget{FUNCTION:Func_GETTRJ}{\texttt{TRJ = Func\_GetTRJ(BD)}}
    \begin{enumerate}[label=(\roman*)]
    \item \verb|DESCRIPTIONs|: \\
    This function allows the user to retrieve the trajectory structure from the bifurcation diagram \verb|BD|.
    \item \verb|INPUTs|:
    \begin{enumerate}[label=(\alph*)]
    \item \verb|BD|: This is a bifurcation diagram structure. Generally, it is stored as a field of an AutoRepo structure. It must contain the field \verb|LABPTs|, and at least one single labeled point must contain a \verb|PO| field. Otherwise, the function fails. For more information regarding the structure of bifurcation diagrams, refer to Sec.~\ref{Section:BifDiag}. 
    \end{enumerate}
    \item \verb|OUTPUTs|: \\
    The output of this function is a structure with multiple trajectories as illustrated in Sec.~\ref{Section:Trajectories}.
    \end{enumerate}

\vspace{0.7cm}

%%%% FUNC_READAUTOREPO()
\noindent \hypertarget{FUNCTION:FUNC_READAUTOREPO}{\texttt{AR = Func\_ReadAutoRepo(M,fAR,"BIFURCATIONs","TRAJECTORIEs")}}
    \begin{enumerate}[label=(\roman*)]
    \item \verb|DESCRIPTION|: \\
    This function allows the user to load the results stored in `.auto' files. It returns an AutoRepo structure with the organization presented in Sec.~\ref{Section:AutoRepo}.
    \item \verb|INPUTs|:
    \begin{enumerate}[label=(\alph*)]
    \item \verb|M|: Structure of the model used to compute the `.auto' file.
    \item \verb|fAR|: Relative or absolute path to the `.auto' file exported from XPPAUT.
    \item \verb|BIFURCATIONs|: This is an optional field, specified via a name-value declaration: \\
    \verb|(...,`BIFURCATIONs',VAL)|, where \verb|VAL| is a boolean value. When it is not specified, its default value is \verb|true|. It is a boolean parameter controlling the loading of the points under the field \verb|LABPTs| of a bifurcation diagram. When this parameter is set to \verb|false| no labeled points and consequently no special trajectories nor periodic orbits are loaded under the \verb|LABPTs| field.
    \item \verb|TRAJECTORIEs|: This is an optional boolean parameter and, similar to \verb|BIFURCATIONs|, it must be passed as a name-value argument. Its default value is \verb|true| and it controls the loading of the \verb|PO| field in the label points. When it is set to \verb|false| the \verb|PO| field is not created in any of the loaded labeled points.
    \end{enumerate}
    \item \verb|OUTPUTs|: \\
    The output of this function is an AutoRepo structure (\verb|AR|). For more information regarding the organization of the previous structure, please refer to Sec.~\ref{Section:AutoRepo}.
    \item \verb|ADDITIONAL NOTEs|
    \begin{enumerate}[label=(\alph*)]
        \item When the function is called, it visualizes the operations' status over the command window. Depending on the value of the optional inputs, \verb|BIFURCATIONs| and  \verb|TRAJECTORIEs|, the content may change.
        \item In most cases, the user does not want to change the optional input values. However, if the results of AUTO continuations generate large `.auto' files due to numerous computed special trajectories, then the user may want to speed up the loading procedure by disabling the detection of the labeled points or only the loading of special trajectories.
    \end{enumerate}
    \end{enumerate}

\vspace{0.7cm}

%%%% FUNC_READDATA()
\noindent \hypertarget{FUNCTION:Func_ReadData}{\texttt{S = Func\_ReadData(M,f)}}
    \begin{enumerate}[label=(\roman*)]
    \item \verb|DESCRIPTION|: \\
    This function is used to load the simulation data exported in the `.dat' file located at the relative or absolute path \verb|f|. The output of this function is the simulation structure \verb|S|.
    \item \verb|INPUTs|:
    \begin{enumerate}[label=(\alph*)]
        \item \verb|M|: Model structure associated with the `.ode' file used to produce the simulation's data.
        \item \verb|f|: Relative or absolute path of a `.dat' file storing a simulation exported from XPPAUT.
    \end{enumerate}
    \item \verb|OUTPUTs|: \\
    The output of this function is a simulation structure \verb|S|. For more details on this structure, please refer to Sec.~\ref{Section:Simulations}. 
    \end{enumerate}

\vspace{0.7cm}

%%%% FUNC_READMODEL()
\noindent \hypertarget{FUNCTION:Func_ReadModel}{\texttt{M = Func\_ReadModel(f)}}
    \begin{enumerate}[label=(\roman*)]
    \item \verb|DESCRIPTION|: \\
    This function extracts the necessary information from the specified `.ode' file in \verb|f| and provides as output a model structure \verb|M|.
    \item \verb|INPUTs|:
    \begin{enumerate}[label=(\alph*)]
        \item \verb|f|: This is a mandatory input and must be an array of characters specifying the relative or absolute path of where the `.ode' file is located.
    \end{enumerate}
    \item \verb|OUTPUTs|: \\
    The output of this function is the model structure \verb|M|. For more details, please refer to Sec.~\ref{Section:Model}.
    \end{enumerate}

\vspace{0.7cm}

%%%% FUNC_READNULLCLINES()
\noindent \hypertarget{FUNCTION:Func_ReadNullclines}{\texttt{NC = Func\_ReadNullclines(fNC)}}
    \begin{enumerate}[label=(\roman*)]
    \item \verb|DESCRIPTION|: \\
    This function parses the nullcline file located at \verb|fNC| and returns a nullcline structure.
    \item \verb|INPUTs|:
    \begin{enumerate}[label=(\alph*)]
    \item \verb|fNC|: Relative or absolute file to the `.dat' file storing the nullclines exported through XPPAUT. Please remember that the name of the exported file must have the format \verb|[text]_x_y.dat|. Where \verb|text| is a sequence of alpha-numerical characters, \verb|x|/\verb|y| is the variable visualized along the x/y-axis in the XPPAUT window from where the nullclines are calculated and exported, respectively.
    \end{enumerate}
    \item \verb|OUTPUTs|: \\
    The output of this function is a nullcline structure. For more information, please refer to Sec.~\ref{Section:Nullclines}.
    \end{enumerate}

\vspace{0.7cm}

%%%% FUNC_VISUALIZEDIAGRAM()
\noindent \hypertarget{FUNCTION:Func_VisualizeDiagram}{\texttt{Func\_VisualizeDiagram(M,BD,"VAR","BRIND","OPTIONs","TYPE")}}
    \begin{enumerate}[label=(\roman*)]
    \item \verb|DESCRIPTION|: \\
    This function visualizes (over the current figure) the bifurcation diagram \verb|BD|. Specifically, in the absence of any other inputs, the visualization obeys the following rules: if it is passed a one-parameter bifurcation diagram as an input argument, the figure is set up to display along the x/y-axis the main bifurcation parameter/first dynamical variable appearing in the model structure \verb|M|. Instead, if it is passed a two-parameter bifurcation diagram, the figure will be set up to display along the x- and the y-axis the main and the secondary bifurcation parameters. The user can manipulate the variables visualized along the x- and the y-axis through the optional field \verb|"VAR"|, which must be specified as a name-value pair.
    \item \verb|INPUTs|:
    \begin{enumerate}[label=(\alph*)]
    \item \verb|M|: Model structure.
    \item \verb|BD|: Bifurcation diagram structure.
    \item \verb|"VAR"|: This is an optional input argument which must be passed via name-value declaration. It accepts a two- or three-element cell array. Each entry and its position in the cell array defines the parameter or dynamical variable to be visualized along one of the x-/y-/z-axis. Other special accepted entries are \verb|L2|/\verb|T|/\verb|F|, which represents the $\mathcal{L}_2$-norm, the period, and the frequency of a point in the BD, respectively.
    \item \verb|"BRIND"|: The indices of branches to be visualized, in the form of a cell array \verb|{b1,b2,...}|. Must be a subset of branches within \verb|BD|. This is an optional input argument.
    \item \verb|"OPTIONs"|: This is an optional input which must be declared via a name-value method. It accepts an option structure created through the function \hyperlink{FUNCTION:Func_DOBD}{\texttt{Func\_DOBD()}}. If not specified, the default visualization options for bifurcation diagrams will be applied. 
    \item \verb|"TYPE"|: This is an optional input that can be specified through a name-value method. It can be set to "Average"/"Initial"/"Upper"/"Lower"/"Standard" to visualize the average/initial/maximum/minimum/maximum and minimum value associated with the periodic solution branches.
    \end{enumerate}
    \item \verb|OUTPUTs|: -
    \end{enumerate}

\vspace{0.7cm}

%%%% FUNC_VISUALIZEEIG()
\noindent \hypertarget{FUNCTION:Func_VisualizeEig}{\texttt{Func\_VisualizeEig(M,BD,"VAR","BRIND","OPTIONs")}}
    \begin{enumerate}[label=(\roman*)]
    \item \verb|DESCRIPTION|: \\
    The default visualization implemented in this function illustrates how the real and the imaginary components of the eigenvalues or the Floquet multipliers vary as the main bifurcation parameter is changed in 1P or 2P-BD. In the same space, a unit circle cylinder highlights where the bifurcations occur. The visualization can be modified by specifying the optional input parameters to this function.
    
    \item \verb|INPUTs|:
    \begin{enumerate}[label=(\alph*)]
    \item \verb|M|: Model structure.% associated with the bifurcation diagram
    \item \verb|BD|: Bifurcation diagram to be visualized.
    \item \verb|"VAR"|: This argument must be passed as a name-value pair. It accepts as input a continuation parameter used to compute the BD specified as input to this function, dynamical variables, and the fields \verb|EigR|/\verb|EigI|. It is used to set up the x-/y- and z-axis in the diagram. This cell array must have at least one continuation parameter and at least one between the \verb|EigR|/\verb|EigI| field. Possible acceptable cell inputs are the following (or permutation of their input elements): \verb|{P1(2),EigR,EigI}|, \verb|{P1,V,EigR(I)}|, \verb|{P1,P2,EigR(I)}|, \verb|{P1,EigR(I)}|, \verb|{EigR,EigI}|. Where $P1/2$ indicates the main/secondary continuation parameter and \verb|V| a generic dynamical variable. 
    
    \item \verb|"BRIND"|: An optional cell array of branch indices to be visualized. Note that bifurcation points (depicted with coloured dots in \hyperlink{FUNCTION:Func_VisualizeDiagram}{\texttt{Func\_VisualizeDiagram()}} are calculated by AUTO to be associated with \textit{one} specific BD branch. Therefore, removing a branch that ends at a bifurcation point may or may not remove the point as well. Similar to \verb|"VAR"|, this argument must be passed as a name-value pair.
    %\item \verb|"TYPE"|: An optional argument specifying whether the real (\verb|"Re"|, default) or imaginary parts (\verb|"Im"|) are visualized. Similar to \verb|"VAR"|, this argument must be passed as a name-value pair.
    \item \verb|"OPTIONs"|: An optional input for visualization options. See \hyperlink{FUNCTION:Func_VisualizeDiagram}{\texttt{Func\_VisualizeDiagram()}} above. Similar to \verb|"VAR"|, this argument must be passed as a name-value pair.
    \end{enumerate}
    \item \verb|OUTPUTs|: --
    \end{enumerate}

\vspace{0.7cm}

%%%% FUNC_VISUALIZELABPOINTS()
\noindent \hypertarget{FUNCTION:Func_VisualizeLabPoints}{\texttt{Func\_VisualizeLabPoints(M,BD,"VAR","BRIND","PTIND","OPTIONs","TYPE")}}
    \begin{enumerate}[label=(\roman*)]
    \item \verb|DESCRIPTION|: \\
    This function visualizes the labeled points of the bifurcation diagram \verb|BD|. By default, the variables visualized along the x- and the y-axis depend on the bifurcation diagram \verb|BD|. Specifically, if \verb|BD| refers to a one-parameter continuation, the x- and the y-axis depict the main bifurcation parameter and the first dynamical variable stored in the model structure M. Instead, if \verb|BD| is a two-parameter bifurcation diagram, the x- and the y-axis visualize the main and the secondary bifurcation parameters. This function accepts two optional inputs, namely \verb|VAR| and \verb|OPTIONs| to customize the visualization.
    \item \verb|INPUTs|:
    \begin{enumerate}[label=(\alph*)]
    \item \verb|M|: Model structure.
    \item \verb|BD|: Bifurcation diagram structure.
    \item \verb|"VAR"|: This is an optional input argument which must be passed via name-value declaration. It accepts a two- or three-element cell array. Each entry and its position in the cell array defines the parameter or dynamical variable to be visualized along one of the x-/y-/z-axis.
    \item \verb|"BRIND"|: The indices of branches whose labeled points are to be visualized, in the form of a cell array \verb|{b1,b2,...}|. Must be a subset of branches within \verb|BD|. Note that each labeled point is assigned to only one branch by AUTO. This is an optional input argument, and it defaults to all branches of \verb|BD|.
    \item \verb|"PTIND"|: Cell array of point indices to be visualized: \verb|{p1,p2,...}|. This is another way (in addition to \verb|"BRIND"|) to select specific labeled points. This is an optional argument with the default set to all points within \verb|BD|. If both \verb|"BRIND"| and \verb|"PTIND"| are specified, the union of the two fields is taken.
    \item \verb|"OPTIONs"|: This is an optional input which must be declared via a name-value method. It accepts an option structure created through the function \hyperlink{FUNCTION:Func_DOBD}{\texttt{Func\_DOBD()}}. If not specified, the default visualization options for bifurcation diagrams will be applied. 
    \item \verb|"TYPE"|: This is an optional input that can be specified through a name-value method. It can be set to "Average"/"Initial"/"Upper"/"Lower"/"Standard" to visualize the average/initial/maximum/minimum/maximum and minimum value associated with the periodic solution branches.
    \end{enumerate}
    \item \verb|OUTPUTs|: -
    \end{enumerate}

\vspace{0.7cm}

%%%% FUNC_VISUALIZENULLCLINES()
\noindent \hypertarget{FUNCTION:Func_VisualizeNullclines}{\texttt{Func\_VisualizeNullclines(NC,"OPTION")}}
    \begin{enumerate}[label=(\roman*)]
    \item \verb|DESCRIPTION|: \\
    This function allows the user to visualize the nullclines saved in the \verb|NC| structure specified as an input argument to this function. The variables visualized along the x- and the y-axis of this chart are those specified in the file name associated with the \verb|NC| structure; as specified in the \hyperlink{FUNCTION:Func_ReadNullclines}{\texttt{Func\_ReadNullclines()}} function.
    \item \verb|INPUTs|:
    \begin{enumerate}[label=(\alph*)]
    \item \verb|NC|: Nullcline structure.
    \item \verb|"VAR"|: An optional input to control which variable to visualize along the x- and the y-axis. It must be specified as a name-value pair. It accepts as input a cell array with two elements containing the names of the variables to visualize along the x-, and the y-axis.
    \item \verb|"OPTION"|: An optional input for visualization options. See \hyperlink{FUNCTION:Func_VisualizeDiagram}{\texttt{Func\_VisualizeDiagram()}}. This argument must be passed as a name-value pair.
    \end{enumerate}
    \item \verb|OUTPUTs|: -
    \end{enumerate}

\vspace{0.7cm}

%%%% FUNC_WRITEPOINTS()
\noindent \hypertarget{FUNCTION:Func_WritePoints}{\texttt{Func\_WritePoints(M,BD,fWP,`VAR')}}
    \begin{enumerate}[label=(\roman*)]
    \item \verb|DESCRIPTIONs|: \\
    This function allows the user to export the bifurcation diagram \verb|BD| passed as input to this function into a format readable from XPPAUT through the command \verb|Graphic| \verb|stuff| $\rightarrow$ \verb|Freeze| $\rightarrow$ \verb|Bif.Diag|.
    \item \verb|INPUTs|:
    \begin{enumerate}[label=(\alph*)]
    \item \verb|M|: Model structure.
    \item \verb|BD|: This is a bifurcation diagram structure. Generally, it is stored as a field of an AutoRepo structure.
    \item \verb|fWP|: Name of the output file.
    \item \verb|VAR|: This is an optional input that must be passed through a name-value declaration. It accepts a cell array with a number of elements equal to two. The first/second entry specifies the variable or parameter to be visualized along the x-/y-axis. If this field is not specified, the variables exported by default are the main bifurcation parameter and the first dynamical variable implemented in the `.ode' file.
    \end{enumerate}
    \item \verb|OUTPUTs|: \\
    This function produces a `.dat' file (parsable through XPPAUT) whose name is specified in the input argument \verb|fWP|.
    \end{enumerate}

\newpage

\bibliographystyle{abbrv}
\bibliography{ref}

\begin{thebibliography}{10}

\bibitem{Bertram1995}
R.~Bertram, M.~Butte, T.~Kiemel, and A.~Sherman.
\newblock Topological and phenomenological classification of bursting oscillations.
\newblock {\em Bulletin of Mathematical Biology}, 57(3):413--439, 1995.

\bibitem{CK1983}
T.~R. Chay and J.~Keizer.
\newblock Minimal model for membrane oscillations in the pancreatic beta-cell.
\newblock {\em Biophysical Journal}, 42(2):181--189, 1983.

\bibitem{xppMATLAB}
R.~Clewley.
\newblock xpp-matlab, 2019.

\bibitem{PyDSTool}
R.~H. Clewley, W.~E. Sherwood, M.~LaMar, and J.~Guckenheimer.
\newblock Pydstool, a software environment for dynamical systems modeling., 2007.

\bibitem{PlotBD}
M.~De~Pitta.
\newblock Plotbd, 2016.

\bibitem{DKO2008CHAOS}
M.~Desroches, B.~Krauskopf, and O.~H. M.
\newblock Mixed-mode oscillations and slow manifolds in the self-coupled fitzhugh-nagumo system.
\newblock {\em Chaos: An Interdisciplinary Journal of Nonlinear Science}, 18(1):015107, 2008.

\bibitem{DKO2008SIAM}
M.~Desroches, B.~Krauskopf, and H.~M. Osinga.
\newblock The geometry of slow manifolds near a folded node.
\newblock {\em SIAM Journal on Applied Dynamical Systems}, 7(4):1131--1162, 2008.

\bibitem{DKO2010}
M.~Desroches, B.~Krauskopf, and H.~M. Osinga.
\newblock Numerical continuation of canard orbits in slow–fast dynamical systems.
\newblock {\em Nonlinearity}, 23(3):739--765, 2010.

\bibitem{MatCont}
A.~Dhooge, W.~Govaerts, Y.~A. Kuznetsov, H.~G.~E. Meijer, and B.~Sautois.
\newblock New features of the software matcont for bifurcation analysis of dynamical systems.
\newblock {\em Mathematical and Computer Modelling of Dynamical Systems}, 14(2):147--175, 2008.

\bibitem{AUTO07p}
E.~J. Doedel, A.~R. Champneys, T.~F. Fairgrieve, Y.~A. Kuznetsov, B.~E. Oldeman, R.~C. Paffenroth, B.~Sandstede, X.~J. Wang, and C.~Zhang.
\newblock Auto-07p: Continuation and bifurcation software for ordinary differential equations, 2007.

\bibitem{XPPAUT}
B.~Ermentrout.
\newblock {\em Simulating, Analyzing, and Animating Dynamical Systems}.
\newblock Society for Industrial and Applied Mathematics, 2002.

\bibitem{xppToolbox}
P.~Fletcher.
\newblock xpptoolbox, 2023.

\bibitem{HH1952}
A.~L. Hodgkin and A.~F. Huxley.
\newblock A quantitative description of membrane current and its application to conduction and excitation in nerve.
\newblock {\em The Journal of Physiology}, 117(4):500--544, 1952.

\bibitem{JohnJP2024}
S.~John, W.~Barnett, A.~Abdala, D.~Zoccal, J.~Rubin, and Y.~Molkov.
\newblock Exploring the role of the kölliker–fuse nucleus in breathing variability by mathematical modelling.
\newblock {\em The Journal of Physiology}, 602(1):93--112, 2024.

\bibitem{JohnChaos2024}
S.~R. John, R.~S. Phillips, and J.~E. Rubin.
\newblock A novel mechanism for ramping bursts based on slow negative feedback in model respiratory neurons.
\newblock {\em Chaos: An Interdisciplinary Journal of Nonlinear Science}, 34(6):063131, 2024.

\bibitem{Kuznetsov2023}
Y.~A. Kuznetsov.
\newblock {\em Numerical Analysis of Bifurcations}.
\newblock Springer International Publishing, 2023.

\bibitem{Martin2024}
M.~Martin and M.~G. Pedersen.
\newblock Modelling and analysis of camp-induced mixed-mode oscillations in cortical neurons: Critical roles of hcn and m-type potassium channels.
\newblock {\em PLOS Computational Biology}, 20(3):1--24, 2024.

\bibitem{MATLAB}
MathWorks.
\newblock Matlab version: 9.14.0 (r2023a), 2023.

\bibitem{Morris1981}
C.~Morris and H.~Lecar.
\newblock Voltage oscillations in the barnacle giant muscle fiber.
\newblock {\em Biophysical journal}, 35(1):193--213, 1981.

\bibitem{Nechyporenko2024}
K.~Nechyporenko, M.~Voliotis, X.~F. Li, O.~Hollings, D.~Ivanova, J.~J. Walker, K.~T. O’Byrne, and K.~Tsaneva-Atanasova.
\newblock Neuronal network dynamics in the posterodorsal amygdala: shaping reproductive hormone pulsatility.
\newblock {\em Journal of The Royal Society Interface}, 21, 2024.

\bibitem{XPPy}
J.~Nowacki.
\newblock Xppy, 2011.

\bibitem{PyXPP}
M.~Olenik.
\newblock Pyxpp, 2018.

\bibitem{Py_XPPCALL}
I.~Prokin and Y.~Park.
\newblock Py\_xppcall, 2017.

\bibitem{Rinzel1985}
J.~Rinzel.
\newblock Bursting oscillations in an excitable membrane model.
\newblock In {\em Ordinary and Partial Differential Equations}, pages 304--316. Springer Berlin Heidelberg, 1985.

\end{thebibliography}

\end{document}